\documentclass[10pt]{article}

\usepackage{amsfonts}
\usepackage{amsxtra}
\usepackage{amsmath}
\usepackage{a4}
\usepackage{amssymb}
\usepackage{theorem}
\usepackage{epsfig}
\newtheorem{theorem}{Theorem}
\newtheorem{lemma}[theorem]{Lemma}
\newtheorem{corollary}[theorem]{Corollary}
{\theorembodyfont{\rmfamily}}
{\theorembodyfont{\rmfamily}}
\newtheorem{proposition}[theorem]{Proposition}
\numberwithin{theorem}{section}
\numberwithin{equation}{section}
\numberwithin{figure}{section}

\newcommand{\given}{\big\vert}

\newcommand{\Ber}{\text{Ber}}

\newcommand{\intersect}{\cap}

\newcommand{\bu}{\mathbf{u}}
\newcommand{\by}{\mathbf{y}}
\newcommand{\bv}{\mathbf{v}}

\newcommand{\RR}{\mathbb{R}}
\newcommand{\PP}{\mathbb{P}}
\newcommand{\ZZ}{\mathbb{Z}}
\newcommand{\NN}{\mathbb{N}}
\newcommand{\EE}{\mathbb{E}\,}

\newcommand{\cC}{\mathcal{C}}

\newcommand{\cB}{\mathcal{B}}
\newcommand{\cS}{\mathcal{S}}
\newcommand{\cT}{\mathcal{T}}

\renewcommand{\Box}{\square}
\newcommand{\ty}{{\tilde{y}}}

\newcommand{\tC}{{\tilde{C}}}
\newcommand{\tu}{{\tilde{u}}}
\newcommand{\tA}{{\tilde{A}}}
\newcommand{\tT}{{\tilde{T}}}
\newcommand{\tM}{{\tilde{M}}}

\newcommand{\tV}{{\tilde{V}}}
\newcommand{\tS}{{\tilde{S}}}

\newcommand{\Ubar}{{\bar{U}}}

\newcommand{\toi}{\rightarrow \infty}

\newcommand{\bx}{\mathbf{x}}
\newcommand{\bfi}{{\bf i}}

\newcommand{\bfu}{{\bf u}}
\newcommand{\bfv}{{\bf v}}

\newcommand{\be}{\mbox{$\mathbb{E}$}}
\newcommand{\bp}{\mbox{$\mathbb{P}$}}

\newcommand{\cf}{\mbox{$\cal F$}}

\newcommand{\br}{\mbox{$\mathbb R$}}
\newcommand{\bn}{\mbox{$\mathbb N$}}

\renewcommand{\thefootnote}{\fnsymbol{footnote}}

\makeatletter
\@addtoreset{equation}{section}

\makeatother

\def\enpf{
     {  \parfillskip=0pt\hfil {\hbox{$\mbox{$\Box$}$}} \par\bigskip  }
     }

\allowdisplaybreaks

\begin{document}

\title{Heavy tails in last-passage percolation}
\author{
\textbf{$\,\,\,\,\,\,\,\,$Ben Hambly}
\hspace{1cm}
\textbf{James B.~Martin}\\
\textit{University of Oxford}
}
\date{March 21, 2006}
\maketitle

\renewcommand{\thefootnote}{}
\footnotetext[1]{
$\!\!\!\!\!$
}

\begin{abstract}
We consider last-passage percolation models in two dimensions,
in which the underlying weight distribution has
a heavy tail of index $\alpha<2$.
We prove scaling laws and asymptotic distributions,
both for the passage times and for the shape of
optimal paths; these are expressed in 
terms of a family (indexed by $\alpha$) of 
``continuous last-passage percolation'' models in the unit square.
In the extreme case $\alpha=0$ (corresponding to 
a distribution with slowly varying tail)
the asymptotic distribution of the optimal path 
can be represented by a random self-similar measure on $[0,1]$,
whose multifractal spectrum we compute.
By extending the continuous last-passage percolation 
model to $\RR^2$ we obtain
a heavy-tailed analogue of the Airy process,
representing the limit of appropriately scaled vectors
of passage times to different points in the plane. 
We give corresponding results for a directed
percolation problem based on $\alpha$-stable L\'evy processes,
and indicate extensions of the results to higher dimensions.
\end{abstract}

\section{Introduction}
Directed last-passage percolation in two dimensions has received much
attention in recent years. In certain specific cases, for example
where the weights at each site are i.i.d.\ with exponential or
geometric distribution, very precise scaling laws and asymptotic
distributions are now known, both for the passage times and for the
shape of optimal paths (see for example \cite{Johshape,
Johanssontransversal, BDMMZ}).  Such cases are closely related to the
longest increasing subsequence problem, and to Markovian interacting
particle systems such as the totally asymmetric exclusion process;
there are also very close links to random matrix theory, for example
to the behaviour of the largest eigenvalue of a large matrix drawn
from the Gaussian Unitary Ensemble (see for example 
\cite{Neilsurvey} for a survey).

It is believed that the behaviour proved for the exponential and
geometric cases should be \textit{universal}, in that the same scaling
laws and asymptotic distributions should occur in the last-passage
percolation model whose underlying weight distribution is from a much
more general class (certainly including any distribution with an
exponentially decaying tail, and maybe also those with sufficiently
light polynomial tails).  The growth models corresponding to these
last-passage percolation problems should belong to the
\textit{Kardar-Parisi-Zhang (KPZ) universality class} (see for example
\cite{KrugSpohn}).  However, only very limited universality results
have been proved: for example, conditions under which laws of large
numbers for the passage times (or ``shape theorems'') hold, and
asymptotics for passage-times close to the boundary of the quadrant
\cite{JBMshape, BodineauMartin, BaiSui}.

In this paper, we study cases in which the tail of the weight
distribution is sufficiently heavy that such shape theorems fail, and
which certainly fall outside the universality class described
above. Specifically, we assume that the tail of the weight
distribution is regularly varying with index $\alpha<2$.  We describe
a family of ``continuous last-passage percolation'' models (indexed by
$\alpha$), and use them to provide scaling laws and asymptotic
distributions for the discrete models, both for the passage times and
for the shape of the optimal paths. Thus we have a universality
result for these heavy-tailed models as the only information
required to determine the scaling limits is the parameter $\alpha$.

One example of an application where such a heavy-tailed assumption is
very natural is in the use of last-passage models to represent
networks of \textit{queues in tandem}.  The vertex weights in the
percolation models correspond to service times in the queueing
systems, and passage times in the percolation models correspond 
to the total time spent in the queueing system by particular customers; 
see for example \cite{GlyWhi, BBM, blocking}.

In Section \ref{main}, we define the discrete last-passage 
percolation model
precisely; we then describe the continuous last-passage
model and state our main convergence results.  We also derive from the
continuous model a stationary process which can be seen as a
heavy-tailed analogue of the \textit{Airy process} (which was
developed by Pr\"ahofer and Spohn \cite{PraSpo} and Johansson
\cite{JohanssonAiry}), and which gives a process limit for vectors of
passage times to different points, appropriately scaled.

The proofs of the convergence results for passage-times are given in
Section \ref{Tsection}, and those for the optimal paths are given in
Section \ref{pathsection}. The results on the heavy-tailed Airy
process are proved in section \ref{Airysection}.

In Section \ref{greedysection} we explore the case where the tail of
the weight distribution is \textit{slowly varying} (i.e.\ $\alpha=0$).
It is no longer possible to provide a non-degenerate limiting
distribution for the passage times; however, asymptotics for the form
of the optimal paths are still possible, and in fact the 
distribution of the limiting path that 
arises can be described in a particularly simple and
algorithmic way. As the path is increasing it can be thought of
as the distribution of a random measure on $[0,1]$;
this measure is self-similar and we compute its multifractal
spectrum.

The \textit{Brownian directed percolation} model has recently been
much studied in various contexts (see for example \cite{GUEs,
GraTraWid, hmo, OcoYor}).  In Section \ref{stablesection} we discuss a
related model in which Brownian motion is replaced by an
$\alpha$-stable L\'evy process, and we again prove distributional
convergence to the continuous last passage percolation model.

The bulk of the paper concerns the case of two-dimensional
last-passage percolation.  However, almost all of the results extend
easily to $d$ dimensions, $d\geq 3$, and now apply for $\alpha<d$.  We
indicate these extensions in Section \ref{dsection}.

Simulations of trees of optimal paths, of the limiting path for
$\alpha=0$, and of heavy-tailed Airy processes are given in Sections
\ref{main},  \ref{greedysection} and \ref{Airysection}.

\section{Main results}\label{main}
\subsection{Definition of the discrete problem}

Let $F$ be a distribution function.
We will assume that the tail of the distribution $F$
is \textit{regularly varying} with index $\alpha\in(0,2)$; that is,
for all $t>0$,
\[
\frac{1-F(tx)}{1-F(x)}\to t^{-\alpha} \text{ as } x\to\infty.
\]
We will also assume (merely for convenience)
that $F$ is a continuous distribution and $F(0)=0$.

The discrete last-passage percolation model with
underlying weight distribution $F$ is usually defined as follows.

Let $X(i,j), i,j\in\NN$ be i.i.d.\ with common distribution $F$.
The quantity $X(i,j)$ 
represents the \textit{weight} at the site $(i,j)\in\ZZ_+^2$.

For $n\in\NN$, we will define the quantity $T^{(n)}$, the 
\textit {last-passage time}
between $(1,1)$ and $(n,n)$.  Let $\Pi_n$ be the set of
\textit{directed paths} between $(1,1)$ and $(n,n)$.  Each such path
begins at $(1,1)$ and ends at $(n,n)$, and each step consists of
increasing one of the two coordinates by 1.
That is, for any $\pi\in\Pi_n$, we can write $\pi=(v_1, v_2, \dots,
v_{2n})$, where $v_1=(1,1)$, $v_{2n}=(n,n)$, and, for each
$i=1,\dots,2n-1$, $v_{i+1}-v_i$ is either $(1,0)$ or $(0,1)$.

The weight of such a path is the sum of the weights $X(i,j)$
associated with the points $(i,j)$ in the path. Then $T^{(n)}$ is the
maximal weight of a directed path between $(1,1)$ and $(n,n)$; that is:
\begin{equation}
\label{Tndef1}
T^{(n)}=\max_{\pi\in\Pi_n} \sum_{v\in\pi} X(v).
\end{equation}
Note that $T^{(n)}$ depends only on the weights $X(v)$,
$v\in\{1,\dots,n\}^2$. 

\subsection{Continuous model}

We start with an alternative representation of the discrete model.
Let $M^{(n)}_1 \geq M^{(n)}_2 \geq \dots \geq M^{(n)}_{n^2}$ be the
order statistics, written in decreasing order, from an i.i.d.\ sample
of size $n^2$ from the distribution $F$.

Consider the set $\{1/n, 2/n, \dots, (n-1)/n, 1\}^2 \subset [0,1]^2$,
of size $n^2$.  Let the sequence $Y^{(n)}_1, \dots, Y^{(n)}_{n^2}$
consist of a random ordering of the points of this set, chosen
uniformly from the $(n^2)!$ possibilities.
We regard $Y^{(n)}_i$ as the location of the $i$th largest weight
$M^{(n)}_i$ (and we have scaled so that all points lie in the box
$[0,1]^2$).

For two points $y, y'\in[0,1]^2$, we say that $y$ and $y'$ are
\textit{compatible}, and write $y\sim y'$, if $y,y'$ are partially
ordered in that \textit{either} $y\leq
y'$ co-ordinatewise, \textit{or} $y'\leq y$
co-ordinatewise. (Informally, one of $y$ and $y'$ is below and to the
left of the other).
An increasing path will consist of a set of points such that every
pair of points in the set is compatible.  We describe the collection
of increasing paths by the collection $\cC^{(n)}$, which depends on
the points $Y^{(n)}_1, \dots, Y^{(n)}_{n^2}$ alone:
\begin{equation}
\label{Cndef}
\cC^{(n)}=
\cC^{(n)}(Y^{(n)}_1,\dots, Y^{(n)}_{n^2})=
\left\{A\subseteq\{1,\dots,n^2\}
\text{ such that for all } i,j\in A, 
Y^{(n)}_i\sim Y^{(n)}_j
\right\}.
\end{equation}

Now we can give a new definition for $T^{(n)}$, equivalent (in
distribution) to (\ref{Tndef1}):
\begin{equation}\label{Tndef2}
T^{(n)}=\max_{A\in\cC^{(n)}}\sum_{i\in A} M^{(n)}_i.
\end{equation}

We formulate the limiting continuous model by defining the
distribution of a random variable $T$ in an analogous way.

First let $Y_1, Y_2, \ldots$ be an i.i.d.\ sequence, with each $Y_i$
uniformly distributed on the square $[0,1]^2$.
Let $W_1, W_2, \ldots$ be an i.i.d.\ sequence of exponential random
variables with mean 1 (independent of the $(Y_i)$).  Now write, for
each $k\in\NN$, $M_k=(W_1+\ldots+W_k)^{-1/\alpha}$.  (Then with
probability 1, $M_k>M_{k+1}$ for each $k$, and $M_k\to0$ as
$k\to\infty$). The motivation for this definition 
is given by equation (\ref{Mconv}) below.
$M_k$ is the $k$th largest weight, which we imagine positioned
at the point $Y_k\in[0,1]^2$. (The set of locations $Y_k$ 
is of course dense in $[0,1]^2$ with probability 1).

Analogously to (\ref{Cndef}), we
represent the set of increasing paths by the collection $\cC$:
\begin{equation}\label{cCdef}
\cC=
\cC(Y_1, Y_2, \dots)=
\left\{A\subseteq\{1,2,\dots\}
\text{ such that for all } i,j\in A, 
Y_i\sim Y_j
\right\}.
\end{equation}
Then define 
\begin{equation}\label{Tdef}
T=\sup_{A\in\cC}\sum_{i\in A} M_i.
\end{equation}

\noindent\textit{Remark:} 
Note that one could equivalently define $T$ in (\ref{Tdef}) 
as the sup of the weight of \textit{finite} increasing paths $A$,
since either the sup is finite in which case 
the weight of any infinite path can be arbitrarily closely
approximated by that of a finite path, or the sup is infinite 
in which case one can find a finite path with an arbitrarily large
weight. In particular
$T$ can be seen as the supremum of a countable family of measurable
random variables, and so is itself measurable. 
In Section \ref{Airysection} an equivalent construction 
of the continuous last-passage problem 
is given using a Poisson random measure approach 
rather than the sequence of i.i.d.\ uniform positions 
in the unit square described above.

\subsection{Convergence results}
First note that for all $k$,
\[ \left( Y^{(n)}_1, Y^{(n)}_2,\dots Y^{(n)}_k \right) \to
\left(Y_1, Y_2,\dots,Y_k\right) \]
in distribution, as $n\to\infty$.

Now define $a_N=F^{(-1)}\left(1-\frac{1}{N}\right)$.  
(As an example,
if the weight distribution $F$ is Pareto($\alpha$), with
$F(x)=1-x^{-\alpha}$, then $a_{N}=N^{1/\alpha}$.  In general,
$\lim_{N\to\infty} \log a_{N}/\log N=1/\alpha$).

Recall that $M^{(n)}_k$ is 
the $k$th largest value from a sample of size $n^2$ 
from the distribution $F$.
We write $\tM^{(n)}_i=a_{n^2}^{-1} M^{(n)}_i$.
Then from classical extreme value theory we have that, for all $k$,
\begin{equation}\label{Mconv}
 \left( \tM^{(n)}_1, \tM^{(n)}_2,\dots \tM^{(n)}_k \right)
\to \left(M_1, M_2,\dots,M_k\right) 
\end{equation}
in distribution, as $n\to\infty$
(see for example Section 9.4 of \cite{Davidbook}).

In particular, $M^{(n)}_i$ is asymptotically of the order of
$a_{n^2}$, for any $i$.  
(For example, for the Pareto distribution $F(x)=1-x^{-\alpha}$
mentioned above, we have
$a_{n^2}=n^{2/\alpha}$).
Since certainly $T^{(n)}\geq M^{(n)}_1$, we
have that $T^{(n)}$ grows asymptotically at least on the order of
$a_{n^2}$. In fact, we will show that this lower bound gives the right
order of magnitude.

Specifically, let $\tT^{(n)}=a_{n^2}^{-1}T^{(n)}
=\sup_{A\in\cC^{(n)}}\sum_{i\in A} \tM^{(n)}_i$.
Then we will show:
\begin{theorem}
\label{Ttheorem}
The random variable $T$ defined at (\ref{Tdef})
is almost surely finite, and 
$\tT^{(n)}\to T$ in distribution as $n\to\infty$.
\end{theorem}

For comparison, one can consider the case of a lighter tail. 
If $\int_0^\infty [1-F(x)]^{1/2}dx<\infty$
(this condition is very slightly stronger than the
existence of a finite second moment) then
a law of large numbers holds:
$n^{-1}T^{(n)}\to\gamma$ as $n\to\infty$
for some deterministic $\gamma$ \cite{animals}.
If the weights are exponential with mean 1,
then $\gamma=4$ \cite{Rost}, and then in fact 
one has the much finer convergence result that
$n^{-1/3}(T^{(n)}-4n)$ converges in distribution
as $n\to\infty$, to the \textit{GUE Tracy-Widom distribution}
\cite{Johshape}.

\begin{figure}[ht]
\centering
\epsfig{figure=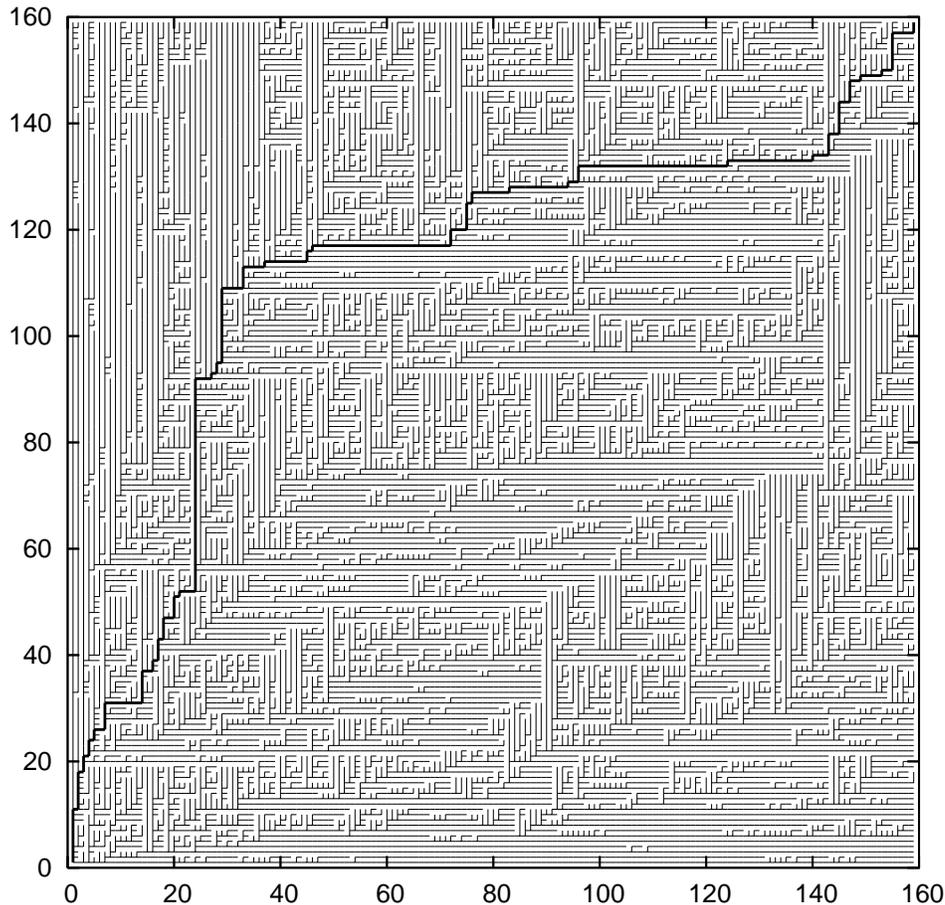, width=1.25\linewidth}
\caption{A simulation of the last-passage percolation
model with $F$ given by a Pareto distribution with index $1$.
The ``tree'' consisting of optimal paths from $(1,1)$
to $(i,j)$, for all $1\leq i,j\leq 159$ is displayed;
the thickened path is the optimal path from $(1,1)$
to $(159, 159)$. This represents $P^{(n)*}$ in the language 
of Theorem \ref{paththeorem2}.
\label{paretofig}}
\end{figure}
\begin{figure}[ht]
\epsfig{figure=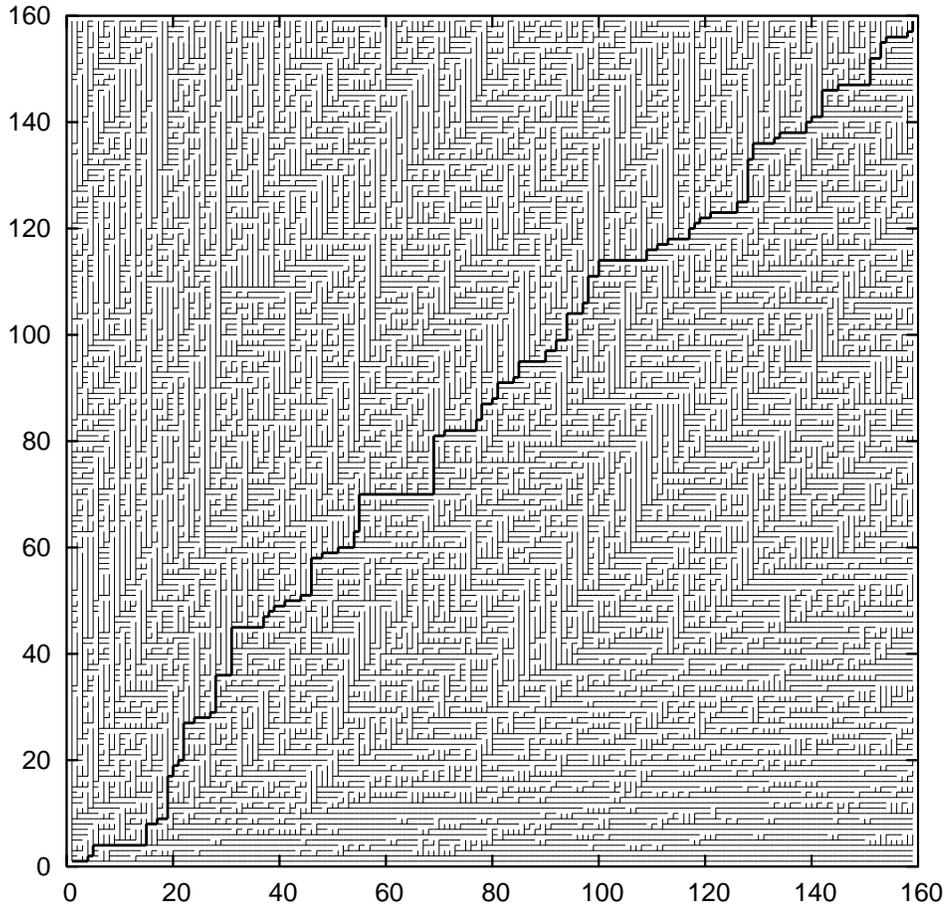, width=1.25\linewidth}
\caption{As Figure \ref{paretofig}, but now with $F$
given by an exponential distribution. As $n$ grows,
the distribution of the optimal path from $(1,1)$ to $(n,n)$ 
becomes concentrated around the diagonal -- the deviations
of the path from a straight line are on the order of $n^{2/3}$.
This contrasts with the deviations on the order of $n$
observed in the heavy-tailed case, where the limiting 
path distribution is non-degenerate (as illustrated in Figure 
\ref{paretofig}).
\label{exponentialfig}}
\end{figure}

We now outline the results on 
path convergence, which are given in full in Section 
\ref{pathsection}.
We will show that the optimal 
path for the continuous model is well defined; 
that is, with probability 1 there exists a unique $A^*\in\cC$
such that $T=\sum_{i\in A^*} M_i$
(attaining the $\sup$ in (\ref{Tdef})). 
Then there is a unique closed connected set $P^*\subset[0,1]^2$
which contains all the points $Y_i, i\in A^*$ 
and which itself has the directed path property 
(i.e.\ $y\sim y'$ for all $y, y'\in P^*$).

Analogously, let $A^{(n)*}$ be the optimal
path for the discrete model, achieving
the $\max$ in (\ref{Tndef2})
(Since the weight distribution $F$
is continuous, the finitely many increasing paths 
all have different weights a.s., and
so this maximizing path is a.s.\ unique). 
Let $P^{(n)*}$ be obtained by linear 
interpolation between the locations 
$\{Y_i^{(n)}: i\in A^{(n)*}\}$
of weights used in this optimal path, taken in increasing order.
Then we show that the distribution of $P^{(n)*}$ converges
to that of $P^*$ as $n\to\infty$ (under the Hausdorff metric
on closed subsets of $[0,1]^2$).

In the corresponding situation for weights with
exponential distribution, the optimal path converges instead 
to a trivial limit, the straight line from $(0,0)$ to $(1,1)$. 
In general, the deviations of the optimal path from 
$(1,1)$ to $(n,n)$ in the discrete model
are expected to be of the order of $n^{2/3}$ in cases
falling into the KPZ universality class 
(proved rigorously in certain cases \cite{Johanssontransversal, BDMMZ}),
rather than on the order of $n$ as we see in the heavy-tailed case.
See Figures \ref{paretofig} and \ref{exponentialfig}
for simulations of optimal paths in the cases
of weight distributions which are Pareto and exponential.

\subsection{Last-passage random fields and the heavy-tailed Airy process}
\label{preAiry}
In Section \ref{Airysection} 
we show how the results described above 
can be extended to give the multivariate convergence 
of vectors of passage times to different points.

For $x, y>0$, let
$T^{(n)}(x,y)$ be the maximal weight of a path from $(1,1)$ to
$(\lceil nx\rceil, \lceil ny\rceil)$.
(So for example the quantity $T^{(n)}$ defined in 
(\ref{Tndef2}) is equal to $T^{(n)}(1,1)$).

Define also $\tT^{(n)}(x,y)=a_{n^2}^{-1}T^{(n)}(x,y)$ as before.

Then we will construct a random field $\{T(x,y), x, y>0\}$, using
a Poisson random measure construction rather than the sequence of
points ordered in decreasing order of weight above, in such a way that
\[ \Big\{\tT^{(n)}(x,y), x, y>0\Big\} \to \Big\{T(x,y),
x, y>0\Big\} \]
as $n\to\infty$, in the sense of convergence of finite-dimensional
distributions. 

From a scaling property of the distribution of the weights
one has further that the random field defined by 
\[ \Theta(u,v)=\exp\left(-\frac{u+v}{\alpha}\right) 
T\left(e^u, e^v\right) \]
is \textit{stationary} on $\RR^2$.
The convergence above can be rewritten as
\begin{equation}
\label{Thetaconv}
\Bigg\{ \exp\left(-\frac{u+v}{\alpha}\right) 
\tT^{(n)}\left(e^u,e^v\right), u,v\in\RR\Bigg\}
\to \Big\{\Theta(u,v), u, v\in\RR\Big\}. 
\end{equation}

\begin{figure}[ht]
\centering
\epsfig{figure=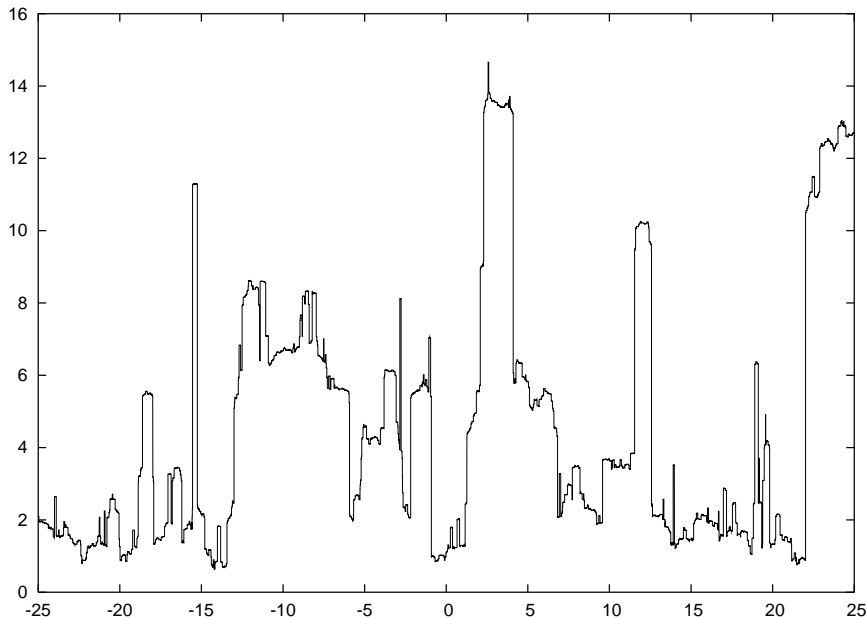, width=0.85\linewidth}
\caption{Figures \ref{Airyfig1}-\ref{Airyfig3}
show simulations of the heavy-tailed Airy process $H_t$ for
three different values of $\alpha$. Here $\alpha=1$.
\label{Airyfig1}}
\end{figure}

\begin{figure}[ht]
\centering
\epsfig{figure=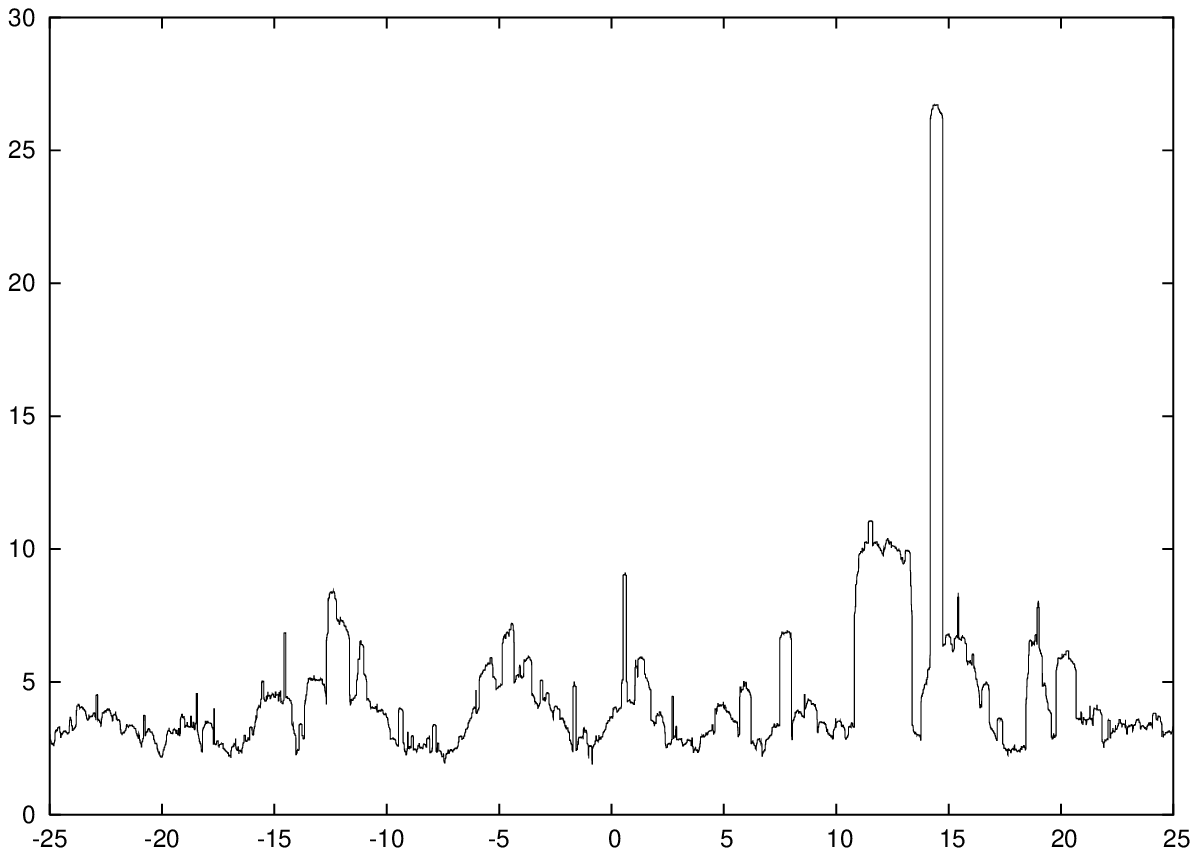, width=0.85\linewidth}
\caption{The heavy-tailed 
Airy process $H_t$ in the case $\alpha=1.5$.\label{Airyfig2}}
\end{figure}

\begin{figure}[ht]
\centering
\epsfig{figure=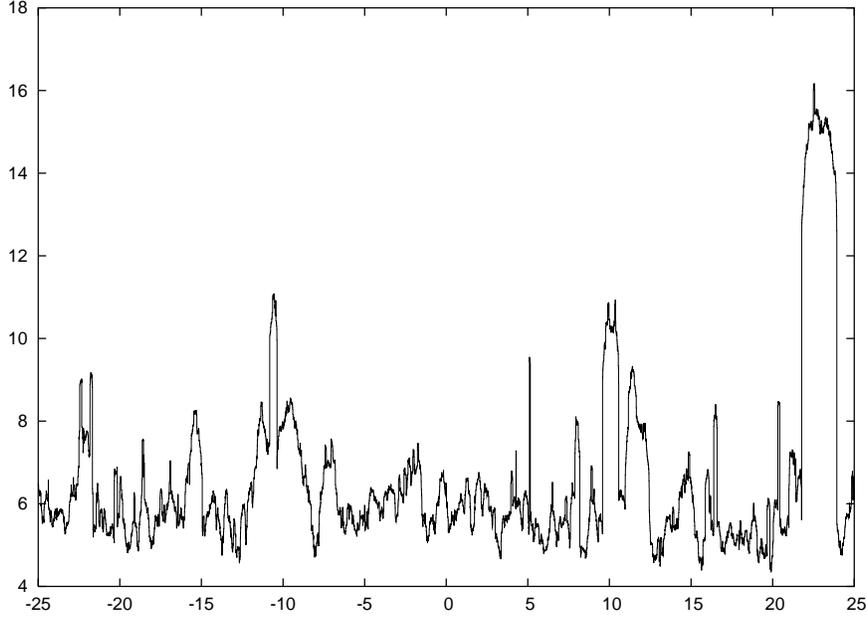, width=0.85\linewidth}
\caption{The heavy-tailed Airy process $H_t$ 
in the case $\alpha=1.99$.\label{Airyfig3}}
\end{figure}

To remove the multiplicative factor
on the LHS of (\ref{Thetaconv}), one can look at a line
$u+v=\text{const}$; for example, the process $H_y=\Theta(y,-y)$.  
This process is stationary in $y\in\RR$, and we obtain the
weak convergence
\[ \left\{ a_{n^2}^{-1} T^{(n)}\left(e^y,e^{-y}\right), y\in\RR\right\} 
\to 
\big\{\Theta(y,-y), y\in\RR \big\}. 
\]
This gives an analogy with the ``Airy process'' \cite{JohanssonAiry}
\cite{PraSpo}, which arises for example in the case where the 
underlying weight
distribution is exponential (with mean 1, say).  
There one obtains 
a stationary process limit
for the quantities
\[
\Big\{n^{-1/3} \Big[ T^{(n)}(1+yn^{-1/3}, 1-yn^{-1/3}) - 4n \Big],
\, y\in\RR \Big\},
\]
whose marginals are given by the GUE Tracy-Widom distribution.
Simulations of the ``heavy-tailed Airy process'' $H_y$ 
are given in Figures \ref{Airyfig1}-\ref{Airyfig3}.

We also obtain estimates on the moments and correlations
of the random field $T$, showing for example that
$\EE T(x,y)^\beta<\infty$ for all $\beta<\alpha$ 
and giving bounds for $\EE|T(x,y)-T(x',y')|^\beta$.

We note that the path convergence described above could
also be extended to the multivariate setting,
to describe the convergence of the distribution of
trees of optimal paths to a continuous tree structure
given by the set of optimal paths in the continuous
last-passage percolation model. However we do not
pursue this further in this paper.

\section{Convergence of the last-passage 
time distribution}
\label{Tsection}

To establish the convergence in Theorem \ref{Ttheorem},
we will work with approximations to $T$ and $T^{(n)}$
which depend only on the $k$ largest weights.
First, define
\begin{gather*}
\cC_k= \cC(Y_1, Y_2, \dots, Y_k)= \left\{A\subseteq\{1,2,\dots,k\}
\text{ such that for all } i,j\in A, Y_i\sim Y_j \right\}. \\
\cC^{(n)}_k= \cC^{(n)}_k(Y^{(n)}_1,\dots, Y^{(n)}_{k\wedge n^2})=
\left\{A\subseteq\{1,\dots,k\wedge n^2\}
\text{ such that for all } i,j\in A, 
Y^{(n)}_i\sim Y^{(n)}_j
\right\}.
\end{gather*}

Note that in fact 
$\cC_k
=\{A\in\cC: A\subseteq\{1,\dots,k\}\}
=\{A\cap\{1,\dots,k\}: A\in\cC\}$, and similarly for $\cC_k^{(n)}$.

Now let
\[
T_k=\sup_{A\in\cC}\sum_{i\in A, i\leq k} M_i,
\]
and 
\[
T^{(n)}_k=\sup_{A\in\cC^{(n)}}\sum_{i\in A, i\leq k} M^{(n)}_i.
\]
Note that indeed $T_k$ depends only on $(M_1,\dots, M_k)$
and $(Y_1,\dots, Y_k)$, while
$T^{(n)}_k$ depends only on $(M^{(n)}_1, \dots, M^{(n)}_k)$
and $(Y^{(n)}_1, \dots, Y^{(n)}_k)$.

As before, define 
$\tT^{(n)}_k=a_{n^2}^{-1}T^{(n)}_k$.

We also define the ``remainder terms'' $S_k$, $S^{(n)}_k$ by
\[
S_k=\sup_{A\in\cC}\sum_{i\in A, i>k} M_i
\,\,
\text{ and }
\,\,
S^{(n)}_k=\sup_{A\in\cC^{(n)}}\sum_{i\in A, i>k} M^{(n)}_i.
\]
Write also $\tS^{(n)}_k=a_{n^2}^{-1}S^{(n)}_k$.

\begin{lemma}
\label{Slemma}
With probability 1,
$S_k<\infty$ for all 
$k\geq 0$, and $S_k\to 0$ as $k\to\infty$.
\end{lemma}

In particular, putting $k=0$ we will have that $T<\infty$ a.s.
(Later on, we will show more, namely that $\EE T^\beta<\infty$
for all $0<\beta<\alpha$; see Proposition \ref{prop:fieldmoment}).

We will also have that $T_k\to T$ a.s.\ as $k\to\infty$,
since, for all $k$,
\begin{align*}
0\leq T-T_k &=\sup_{A\in\cC} \sum_{i\in A} M_i - \sup_{A\in\cC} \sum_{i\in A,
  i\leq k} M_i \\
&\leq \sup_{A\in\cC} \sum_{i\in A, i>k} M_i \\
&= S_k.
\end{align*}

The convergence in Theorem \ref{Ttheorem} will then follow from the
following two results, which provide control over $T_k - \tT^{(n)}_k$
and $\tT^{(n)}_k - \tT^{(n)}$ for appropriate $k$:

\begin{proposition}
\label{coupledprop}
Let $\epsilon>0$ and $k$ be fixed.  Then for all 
$n$ sufficiently large,
say $n\geq N_k(\epsilon)$, there is a coupling 
of the continuous model and the discrete model 
indexed by $n$ under which 
\begin{gather}
\label{coupled1}
\PP\left(\sum_{i=1}^k \left| M_i - \tM^{(n)}_i \right| > \epsilon
\right) \leq \epsilon, \\
\label{coupled2}
\PP\left( \sum_{i=1}^k \left\| Y_i - Y^{(n)}_i \right\| > \epsilon
\right) \leq \epsilon, \\
\label{coupled3}
\PP\left( \cC_k^{(n)}\neq \cC_k \right) \leq \epsilon.
\end{gather}
\end{proposition}

\begin{proposition}
\label{Aprop}
Let $\epsilon>0$. Then for $k$ sufficiently large,
\[
\PP\left(\tS_k^{(n)}>\epsilon\right)\leq \epsilon
\]
for all $n$.
\end{proposition}

\noindent\textbf{Proof of Lemma \ref{Slemma}}:

First, we define $L_i=\sup_{A\in\cC}\left|A\cap\{1,\dots,i\}\right|$.
$L_i$ is the largest number of the points $Y_1,\dots,Y_i$ (the
locations of the $i$ largest weights) that can be included in an
increasing path.  Note that the collection $(L_i)$ is independent of
the collection $(M_i)$. 
$L_i$ has the distribution of the ``longest increasing subsequence''
of a random permutation of length $i$. In particular, there is a 
constant $c$ such that, for all $i$, $\EE L_i \leq c\sqrt{i}$ 
and $\EE L_i^2 \leq c i$; also, 
$L_i/\sqrt{i}\to 2$ in distribution.
See for example \cite{AldDiapatience} for a survey.

We will also write $U_k=\sum_{i=k+1}^\infty L_i(M_i-M_{i+1})$ for each
$k\geq 0$.  Fix $A\in\cC$, and define $R_i=|A\cap\{1,\dots,i\}|$.
Then $I(i\in A)=R_i-R_{i-1}$, and by definition $R_i\leq L_i$.

We have
\begin{align*}
\sum_{i\in A, i>k} M_i
&=\lim_{n\to\infty}\sum_{i\in A, k<i\leq n} M_i
\\
&=\lim_{n\to\infty}\sum_{i=k+1}^n M_i I(i\in A)
\\
&=\lim_{n\to\infty}\sum_{i=k+1}^n M_i(R_i-R_{i-1})
\\
&=\lim_{n\to\infty}
\left[
-M_{k+1}R_k + \sum_{i=k+1}^{n-1} R_i(M_i-M_{i+1})
+M_n R_n
\right]
\\
&\leq
\lim_{n\to\infty}
\sum_{i=k+1}^{n-1} R_i(M_i-M_{i+1}) 
+\liminf_{n\to\infty} M_n R_n
\\
&\leq
\lim_{n\to\infty}
\sum_{i=k+1}^{n-1} L_i(M_i-M_{i+1}) 
+\liminf_{n\to\infty} M_n L_n
\\
&=
U_k +\liminf_{n\to\infty} M_n L_n.
\end{align*}


Now 
$\liminf_{n\to\infty} M_n L_n = 0$ a.s.;
this follows, for example, since 
(by the law of large numbers) $M_n \sim n^{-1/\alpha}$ a.s.\
(with $\alpha<2$), 
and since $L_n/\sqrt{n}$ converges in distribution to a constant.

Since the inequality above holds for any $A\in\cC$, 
we therefore have that $S_k\leq U_k$ a.s., for any $k$.
To conclude the proof, we will show that with probability 1,
$U_k<\infty$ for all $k$, and $U_k\to 0$ as $k\to\infty$.
(In fact, as soon as $U_1<\infty$, we necessarily 
have that $U_k\to0$ as $k\to\infty$,
since the quantity $U_k$ is the ``remainder'' from 
index $k+1$ onwards in the infinite sum $U_1$; 
if the infinite sum $U_1$ converges, then by definition
these remainders tend to 0).

Hence it's enough that $U_k<\infty$ a.s., for all $k$.
Specifically, we'll show that 
$\EE U_k$ is finite whenever $k>1/\alpha$.
Then certainly $U_k<\infty$ a.s.\ for such $k$, and
in fact $U_r<\infty$ a.s.\ for all $r$, since if $r<k$,
$U_r-U_k$ is the sum of only finitely many terms.

By independence of the collections $(L_i)$ and $(M_i)$,
\begin{align}
\nonumber
\EE U_k
&=\sum_{i=k+1}^\infty \EE L_i(\EE M_i - \EE M_{i+1})
\\
\label{Sbound}
&\leq 
\sum_{i=k+1}^\infty ci^{1/2}(\EE M_i - \EE M_{i+1}).
\end{align}

Now $M_r$ has the distribution of $(V_r)^{-1/\alpha}$, 
where $V_r$ has Gamma($r,1$) distribution. We then obtain
\begin{align}
\EE M_r
\nonumber
&=\int_0^\infty \frac{1}{\Gamma(r)} v^{r-1}e^{-v}v^{-1/\alpha} dv
\\
\label{Gammacalc}
&=\Gamma\big(r-1/\alpha\big)/\Gamma(r).
\end{align}
Using the identity $\Gamma(z+1)=z\Gamma(z)$ and the
fact that the gamma function is log convex,
one has
\begin{equation}
\label{Gammafact}
(x-1)^a\leq \frac{\Gamma(x+a)}{\Gamma(x)} \leq (x+a)^a
\end{equation}
for $x>1$, $a<0$.
Then 
\begin{align*}
\EE M_r-\EE M_{r+1}
&=\frac{\Gamma(r-1/\alpha)}{\Gamma(r)}-\frac{\Gamma(r+1-1/\alpha)}{\Gamma(r+1)}
\\
&=\frac{\Gamma(r-1/\alpha)\left[1-\frac{r-1/\alpha}{r}\right]}{\Gamma(r)}
\\
&=\frac{1}{\alpha r}\frac{\Gamma(r-1/\alpha)}{\Gamma(r)}
\\
&\leq \frac{1}{\alpha r}\big(r-1/\alpha-1)^{-1/\alpha}.
\end{align*}

Returning to (\ref{Sbound}), we have
\[
\EE U_k\leq \frac{c}{\alpha}\sum_{i=k+1}^\infty i^{-1/2}
\big(i-1/\alpha-1)^{-1/\alpha},
\]
which is finite for all $k>1/\alpha$ (since $\alpha<2$).$\hfill\Box$

\medskip

\noindent\textbf{Proof of Theorem \ref{Ttheorem}}:
We will find a coupling of $\tilde{T}^{(n)}$ and $T$ for each $n$
such that $\tilde{T}^{(n)}-T\to 0$ in probability as $n\to\infty$.

For each $n\in\NN$, define $k_n=\max\{k:n\geq N_k(1/k)\}$.

Then $k_n\to\infty$ as $n\to\infty$, and, for all $n$,
$n\geq N_{k_n}(1/k_n)$.
Hence from Propositions \ref{coupledprop} and \ref{Aprop}
and from Lemma \ref{Slemma},
there are couplings such that, as $n\to\infty$,
\begin{gather}
\label{conv}
\sum_{i=1}^{k_n} \left| M_i - \tM^{(n)}_i \right| \to 0
\\
\nonumber
\sum_{i=1}^{k_n} \left\| Y_i - Y^{(n)}_i \right\| \to 0
\\
\nonumber
\tS^{(n)}_{k_n}\to 0
\\
\nonumber
S_{k_n}\to 0
\\
\intertext{in probability, and}
\nonumber
\PP\left( \cC_{k_n}^{(n)}\neq \cC_{k_n} \right) \to 0.
\end{gather}

Now 
\[ T-\tT^{(n)} = \left( T-T_{k_n} \right) + \left( T_{k_n} -
\tT^{(n)}_{k_n} \right) + \left( \tT^{(n)}_{k_n} - \tT^{(n)}. \right)
\] 
We have $|T-T_{k_n}|\leq S_{k_n}$ and $|\tT^{(n)}_{k_n}-\tT^{(n)}|\leq
\tS^{(n)}_{k_n}$, so to show that the LHS converges to 0 in
probability as desired, it remains to show that $\left(
T_{k_n}-\tT^{(n)}_{k_n} \right)\to 0$ in probability.

We have
\[
T_{k_n}=\max_{A\in\cC_{k_n}} \sum_{i\in A} M_i
\qquad
\text{ and }
\qquad
\tT^{(n)}_{k_n}=\max_{A\in\cC^{(n)}_{k_n}} \sum_{i\in A} \tM^{(n)}_i,
\]
so if $\cC_{k_n}=\cC^{(n)}_{k_n}$, then
\[
\left| T_{k_n}-\tT^{(n)}_{k_n} \right|
\leq
\sum_{i=1}^{k_n}\left| M_i-\tM^{(n)}_i \right|.
\]
Since $\PP\left(\cC_{k_n}\neq\cC^{(n)}_{k_n}\right)\to 0$ 
and $\sum_{i=1}^{k_n}\left| M_i-\tM^{(n)}_i \right|\to 0$
in probability, we are done.\enpf

To complete the proof it remains to prove
Propositions \ref{coupledprop} and \ref{Aprop}.

\subsection{Convergence of $\tT^{(n)}_k$ to $T_k$}

\noindent\textbf{Proof of Proposition \ref{coupledprop}}:

We have 
$(\tM^{(n)}_1,\dots,\tM^{(n)}_k,Y^{(n)}_1,\dots,Y^{(n)}_k)
\to
(M_1,\dots,M_k,Y_1,\dots,Y_k)$ in distribution as $n\to\infty$.
By the Skorohod Representation Theorem, we can define 
all the variables on the same space in such a way that
the convergence occurs almost surely.
Then indeed (\ref{coupled1}) and (\ref{coupled2}) must hold for 
large enough $n$.

Note that since the variables $Y_i$ are i.i.d.\ uniform 
on $[0,1]^2$, there are almost surely no two $i$ and $j$
such that $Y_i(d)=Y_j(d)$ for $d=1$ or $2$.

Thus if we perturb the point $(Y_1,\dots, Y_k)$ by a small 
enough amount, the orderings of all the coordinates 
remain the same, and the set $\cC_k$ of increasing paths is unchanged.
In fact, if
\[
\max_{1\leq i\leq k} \|Y_i-Y^{(n)}_i\|
\leq
\frac{1}{2}\min_{1\leq i,j \leq k, i\neq j} \min_{d=1,2} |Y_i(d)-Y_j(d)|,
\]
then $\cC^{(n)}_k=\cC_k$. Since we have $Y^{(n)}_i\to Y_i$
a.s.\ on the joint probability space for all $1\leq i\leq k$,
we then have $\cC^{(n)}_k=\cC_k$ eventually,
with probability 1.
Thus (\ref{coupled3}) must also hold for all large enough $n$,
as desired. $\hfill\Box$

\subsection{Convergence of $\tT^{(n)}-\tT^{(n)}_k$ to $0$}

Our aim in this section is to prove Proposition \ref{Aprop}.

Define the ``good event'' $\cB^{(n)}_k$:
\begin{equation}
\label{Bdef}
\cB^{(n)}_k =\left\{
F^{-1}\left(1-\frac{2r}{n^2}\right)
\leq
M^{(n)}_r
\leq 
F^{-1}\left(1-\frac{1}{n^2}\right)
\text{ for all } k<r\leq n^2
\right\}.
\end{equation}

\begin{lemma}
\label{Blemma}
$\PP\left(\cB^{(n)}_k\right)\to 1$ as $k\to\infty$, uniformly in $n$.
\end{lemma}

\noindent\textit{Proof:}
\begin{align*}
\PP\left(\cB^{(n)}_k \text{ fails }\right)
&\leq \PP\left( M^{(n)}_{k+1}> F^{-1}\left(1-\frac{1}{n^2}\right)
\right) +\sum_{r=k+1}^{\lfloor n^2/2 \rfloor}
\PP\left( M^{(n)}_r< F^{-1}\left(1-\frac{2r}{n^2}\right) \right)
\\
&= \PP\left(\text{Binomial}(n^2,1/n^2)\geq k+1\right)
+\sum_{r=k+1}^{\lfloor n^2/2 \rfloor}
\PP\left(\text{Binomial}(n^2, 2r/n^2) < r\right)
\\
&\leq\frac{1}{k+1} + \sum_{r=k+1}^{\lfloor n^2/2 \rfloor} 
2\exp\left(-\frac{r}{12}\right),
\end{align*}
using Markov's inequality for the first term and 
an estimate from Corollary 2.3 of \cite{JLRbook} for the second.
The RHS tends to 0 as $k\to\infty$, uniformly in $n$, as 
required.$\hfill\Box$

\medskip

Now, we will prove a bound on the expectation
of $\tS^{(n)}_k$ 
in terms of the ``order statistics''
$M^{(n)}_r$. 
Define
\[
L^{(n)}_i= \max_{A\in\cC^{(n)}} \left| A\cap\{1,2,\dots,i\} \right|.
\]
Recall that $Y^{(n)}_r\in\{1,2,\dots,n\}^2$ is
the location of the $r$th largest weight, $M^{(n)}_r$.
Thus, $L^{(n)}_i$ is the maximum number of the 
points $Y^{(n)}_1,\dots,Y^{(n)}_i$ that can be included
in an increasing path.


Note that the collection $(L^{(n)}_r)_{1\leq r\leq n^2}$ 
is a function of the values $Y^{(n)}_r$ alone;
in particular it is independent of the weights $M^{(n)}_r$
and of the events $B^{(n)}_r$.

\begin{lemma}
\label{Lboundlemma}
There is a constant $c$ independent of $m$ and $n$ such 
that $\EE L^{(n)}_m \leq c\sqrt{m}$, whenever $1\leq m\leq n^2$.
\end{lemma}

\noindent\textit{Proof:}
The distribution of 
$\{Y^{(n)}_1,\dots,Y^{(n)}_m\}$
is uniform over the subsets of $\{1,\dots,n\}^2$ of size $m$,
and 
$L^{(n)}_r$ is the maximum number of the 
points $Y^{(n)}_1,\dots,Y^{(n)}_r$ that can be included
in an increasing path.

We compare this with the last-passage percolation
problem in $\{1,\dots,n\}^2$ with i.i.d.\
Bernoulli($p$) weights.

We have the following representation for the expectation 
of the passage time for such a problem:
\begin{equation}
\label{dec}
\EE_{\Ber(p)} T(n,n)=\sum_{r} \PP(B_{p,n^2}=r)\EE L^{(n)}_r,
\end{equation}
where 
$B_{p,n^2}$ has Binomial($p,n^2$) distribution,
since, conditional on the event
that exactly $r$ of the $n^2$ weights
have value 1, the set of positions of the weights 
with value 1 is uniformly distributed among all the 
subsets of $\{1,\dots,n\}^2$ of size $r$.

Proposition 2.2 of \cite{animals} shows that there
exists a constant $c_1$ such that 
\begin{equation}
\label{animalsbound}
\EE_{\Ber(p)} T(n,n)\leq c_1 p^{1/2}n
\end{equation}
for all $n$, $p$.

Given $1\leq m\leq n^2$, set $p=\min(2m/n^2,1)$.
Using (\ref{dec}), (\ref{animalsbound}) and the 
fact that $\EE L^{(n)}_r$ is increasing in $r$,
we obtain
\begin{align*}
\PP(B_{p,n^2}\geq m)\EE L^{(n)}_m
&\leq \sum_{r} \PP(B_{p,n^2}=r)\EE L^{(n)}_r
\\
&\leq c_1 \left(\min(2m/n^2,1)\right)^{1/2} n
\\
&\leq c_2\sqrt{m}.
\end{align*}

To complete the proof, it then suffices to bound
$\PP(B_{p,n^2}\geq m)$ away from 0 uniformly in $1\leq m\leq n^2$.

If $m\geq n^2/2$ then $p=1$ and $\PP(B_{p,n^2}\geq m)=1$.

For $1\leq m\leq n^2/2$, we have $pn^2/2=m$,
and we use the estimate
\[
\PP(B_{p,n^2}<pn^2/2)\leq \exp(-pn^2/8),
\]
(see for example Theorem 2.1 of \cite{JLRbook}),
to give
$\PP(B_{p,n^2}\geq m)\geq 1-\exp(-1/8)$ uniformly 
in $1\leq m\leq n^2/2$ 
as desired.$\hfill\Box$

\begin{lemma}
\label{differencelemma}
\[
\EE\left(\tS^{(n)}_k;\cB^{(n)}_k\right)
\leq
c(k+1)^{1/2}\EE \left(\tM^{(n)}_{k+1};\cB^{(n)}_k\right) 
+ c\sum_{r=k+2}^{n^2} r^{-1/2} \EE \left(\tM^{(n)}_r;\cB^{(n)}_k\right).
\]
\end{lemma}

\noindent\textit{Proof:}
The argument is similar to the proof of Lemma \ref{Slemma}.

Let $\tA$ achieve the max in the definition of $\tS^{(n)}_k$,
so that $\tS^{(n)}_k=\sum_{i\in\tA, i>k}\tM^{(n)}_i$.

Define
$R_i=|\tA\cap\{1,2,\dots,i\}|$ for each $i$.
Then $R_i-R_{i-1}=I(i\in\tA)$,
and by definition $R_i\leq L_i^{(n)}$ for each $i$.
We then have
\begin{align*}
\tS^{(n)}_k
&=
\sum_{i\in\tA, i>k} \tM^{(n)}_i
\\
&= \sum_{i=k+1}^{n^2} \tM^{(n)}_i\left(R_i-R_{i-1}\right)
\\
&= -R_k \tM^{(n)}_{k+1}
+\sum_{i=k+1}^{n^2-1} R_i\left(\tM^{(n)}_i-\tM^{(n)}_{i+1}\right)
+R_{n^2} \tM^{(n)}_{n^2}
\\
&\leq 
\sum_{i=k+1}^{n^2-1} L^{(n)}_i\left(\tM^{(n)}_i-\tM^{(n)}_{i+1}\right)
+L^{(n)}_{n^2} \tM^{(n)}_{n^2},
\end{align*}
since $R_i^{(n)}\leq L_i^{(n)}$ and $\tM_i^{(n)}\geq \tM_{i+1}^{(n)}$.

We now take expectations, restricted to the event 
$\cB^{(n)}_k$, using 
Lemma \ref{Lboundlemma} and
the independence of the $L^{(n)}_r$ 
from the $\tM^{(n)}_r$:
\begin{align*}
\EE\left(\tS^{(n)}_k;\cB^{(n)}_k\right)
&\leq
\sum_{i=k+1}^{n^2-1} 
\EE L^{(n)}_i
\left[\EE\left(\tM^{(n)}_i;\cB^{(n)}_k\right)
-\EE\left(\tM^{(n)}_{i+1};\cB^{(n)}_k\right)\right]
+\EE L^{(n)}_{n^2} 
\EE\left(\tM^{(n)}_{n^2};\cB^{(n)}_k\right)
\\
&\leq
\sum_{i=k+1}^{n^2-1} 
c\sqrt{i}
\left[\EE\left(\tM^{(n)}_i;\cB^{(n)}_k\right)
-\EE\left(\tM^{(n)}_{i+1};\cB^{(n)}_k\right)\right]
+
cn
\EE L^{(n)}_{n^2} 
\EE\left(\tM^{(n)}_{n^2};\cB^{(n)}_k\right)
\\
&=
c\sqrt{k+1}\EE \left(\tM^{(n)}_{k+1};\cB^{(n)}_k\right)
+ 
c\sum_{i=k+2}^{n^2} 
\left(\sqrt{i}-\sqrt{i-1}\right)\EE \left(\tM^{(n)}_i;\cB^{(n)}_k\right)
\\
&\leq
c\sqrt{k+1}\EE \left(\tM^{(n)}_{k+1};\cB^{(n)}_k\right)
+ 
c\sum_{i=k+2}^{n^2} 
i^{-1/2}
\EE \left(\tM^{(n)}_i;\cB^{(n)}_k\right),
\end{align*}
since $i^{1/2}-(i-1)^{1/2}\leq i^{-1/2}$. 
This is the required result.$\hfill\Box$

\medskip

The next lemma gives an estimate
on the tail behaviour of the weight distribution
using the regular variation condition.
Note that if the weight distribution were Pareto($\alpha$),
then $F^{-1}(u)=(1-u)^{-1/\alpha}$,
and one then has exactly 
$F^{-1}(u_1)= 
F^{-1}(u_0)\left[
(1-u_1)/(1-u_0)
\right]^{-1/\alpha}$ for any $u_0$, $u_1$.

\begin{lemma}
\label{regboundlemma}
For any $\delta>0$, there exists $U(\delta)<1$ such that
for all $u_0, u_1$ with $U(\delta) \leq u_1 \leq u_0$, 
\[
F^{-1}(u_1)\leq 
2F^{-1}(u_0)\left(
\frac{1-u_1}{1-u_0}
\right)^{-\frac{1}{\alpha}+\delta}.
\]
\end{lemma}

\noindent\textit{Proof:}
Fix $s>1$ sufficiently small that $s^{1/\alpha-\delta}<2$.

From the fact that the tail of $F$ is regularly
varying with index $\alpha$, the following property holds:
if $u_1$ is sufficiently close to $1$ (at least $U(\delta)$, say),
then for all $(1-u_0)<(1-u_1)/s$, 
\[
\frac{F^{-1}(u_1)}{F^{-1}(u_0)} < s^{-\frac{1}{\alpha}+\delta}.
\]

Iterating, one obtains that if $U(\delta)\leq u_1\leq u_0$, then
\begin{align*}
\frac{F^{-1}(u_1)}{F^{-1}(u_0)} 
&< 
\left(
s^{-\frac{1}{\alpha}+\delta}
\right)^{\left\lfloor \log_s \left((1-u_1)/(1-u_0)\right)\right\rfloor}
\\
&=
\exp\left[
\left(-\frac{1}{\alpha}+\delta\right)(\log s)
\genfrac\lfloor\rfloor{}{}{\log\left((1-u_1)/(1-u_0)\right)}{\log s}
\right]
\\
&\leq
\exp\left[
\left(-\frac{1}{\alpha}+\delta\right)
\left(\log\frac{1-u_1}{1-u_0}-\log s\right)
\right]
\\
&=
\genfrac(){}{}{1-u_1}{1-u_0}^{-\frac{1}{\alpha}+\delta}
s^{\frac{1}{\alpha}-\delta}
\\
&\leq
2\genfrac(){}{}{1-u_1}{1-u_0}^{-\frac{1}{\alpha}+\delta}
\end{align*}
as required.$\hfill\Box$

\medskip

Finally we use this tail estimate to control
the expectation of the variables $\tM^{(n)}_r$,
restricted to the ``good set'' $\cB^{(n)}_k$.

\begin{lemma}
\label{EMlemma}
Let $\delta>0$. Then there exist $c_0$, $c_1$ and $c_2>0$ such that
\[
\EE\left(
\tM^{(n)}_r; \cB^{(n)}_k 
\right)
\leq
c_0  r^{-\frac{1}{\alpha}+\delta}
+ c_1 a_{n^2}^{-1} I(r\geq c_2 n^2)
\]
for all $n$, $k$, $r$ satisfying 
$2(1+1/\alpha)< k < r \leq n^2$.
\end{lemma}

\noindent\textit{Proof:}
Let $U(\delta)$ be as in Lemma \ref{regboundlemma},
and set $c_2 = (1-U(\delta))/2$.
Then 
$r<c_2 n^2 \Leftrightarrow 1-2r/n^2 > U(\delta)$.

Note that if 
\[
\max\left\{U(\delta),1-\frac{2r}{n^2}\right\} \leq u \leq 1-\frac{1}{n^2},
\]
then, by Lemma \ref{regboundlemma},
\begin{align}
\nonumber
F^{-1}(u)
&\leq 
2F^{-1}\left(1-\frac{1}{n^2}\right)
\left[n^2(1-u)\right]^{-\frac{1}{\alpha}+\delta}
\\
\nonumber
&= 2 a_{n^2} n^{-2/\alpha}(1-u)^{-1/\alpha}\left[n^2(1-u)\right]^\delta
\\
\label{powerbound}
&\leq
2 a_{n^2} n^{-2/\alpha}(1-u)^{-1/\alpha}(2r)^\delta.
\end{align}

Since $r>k$, we have that
\[
\cB^{(n)}_k\subseteq
\left\{
F^{-1}\left(1-\frac{2r}{n^2}\right)
\leq
M^{(n)}_r
\leq
F^{-1}\left(1-\frac{1}{n^2}\right)
\right\}.
\]

Hence
\begin{align}
\nonumber
\EE\Big(\tM^{(n)}_r ; &\cB^{(n)}_k\Big)
\\
&\leq a_{n^2}^{-1} \EE\left(M^{(n)}_r ; 
F^{-1}\left(1-\frac{2r}{n^2}\right)
\leq M^{(n)}_r \leq F^{-1}\left(1-\frac{1}{n^2}\right)
\right)
\\
\nonumber
&= a_{n^2}^{-1} \EE\left(F^{-1}\left(U^{(n)}_r\right) ; 
1-\frac{2r}{n^2} \leq U^{(n)}_r \leq 1-\frac{1}{n^2} \right)
\\ \nonumber &= a_{n^2}^{-1}
\int_{1-2r/n^2}^{1-1/n^2} F^{-1}(u) f_{r;n^2}(u) du
\\ \nonumber &= a_{n^2}^{-1}
I\left\{1-\frac{2r}{n^2} \leq U(\delta)\right\}F^{-1}(U(\delta))
+a_{n^2}^{-1}
\int_{\max\left\{U(\delta),1-2r/n^2\right\}}^{1-1/n^2} 
F^{-1}(u) f_{r;n^2}(u) du
\\
\nonumber
&\leq
c_1 a_{n^2}^{-1}
I\left(r\geq c_2 n^2\right)
+\int_{\max\left\{U(\delta),1-2r/n^2\right\}}^{1-1/n^2} 
2 n^{-2/\alpha}(1-u)^{-1/\alpha}(2r)^\delta
f_{r;n^2}(u) du
\\
\label{midstep}
&\leq
c_1 I\left(r\geq c_2 n^2\right)
+c_3 a_{n^2} n^{-2/\alpha}r^{\delta}
\int_0^1 (1-u)^{-1/\alpha}f_{r;n^2}(u) du,
\end{align}
where $c_1=F^{-1}(U(\delta))$ and $c_3=2^{1+\delta}$,
and where $f_{r;n^2}$ is the density function 
of the $r$th largest from an i.i.d.\ sample of size $n^2$ from 
the uniform distribution on $[0,1]$.

Now
\[
f_{r;n^2}(u)=
\frac{\Gamma(n^2+1)}{\Gamma(n^2-r+1)\Gamma(r)}
(1-u)^{r-1}u^{n^2-r},
\]
and, since $r-1-1/\alpha>0$, one then has
\begin{align}
\nonumber
\int_0^1 (1-u)^{-1/\alpha}f_{r;n^2}(u) du
&=
\frac{\Gamma(n^2+1)}{\Gamma(n^2-r+1)\Gamma(r)}
\int_0^1 (1-u)^{r-1-1/\alpha}u^{n^2-r} du
\\
\nonumber
&=
\frac{\Gamma(n^2+1)}{\Gamma(n^2-r+1)\Gamma(r)}
\genfrac(){}{}
{\Gamma(n^2+1-1/\alpha)}
{\Gamma(n^2-r+1)\Gamma(r-1/\alpha)}^{-1}
\\
\label{Gammaid}
&=
\frac{\Gamma(n^2+1)}{\Gamma(n^2+1-1/\alpha)}
\frac{\Gamma(r-1/\alpha)}{\Gamma(r)}.
\end{align}

Using (\ref{Gammafact}),
the RHS of (\ref{Gammaid}) is bounded above by
\[
\frac{(n^2+1)^{1/\alpha}}{(r-1-1/\alpha)^{1/\alpha}}.
\]
Since $r>2(1+1/\alpha)$, this is in turn 
no greater than $\left(4n^2/r\right)^{1/\alpha}$.
Inserting this into (\ref{midstep}) gives the desired 
result.$\hfill\Box$

\medskip

\noindent\textbf{Proof of Proposition \ref{Aprop}:}
We may assume that $k\leq n^2$, 
since if $k>n^2$ then $\tS^{(n)}_k=0$.

Fix $\epsilon>0$, and fix some $\delta<\frac{1}{\alpha}-\frac{1}{2}$.
Using Markov's inequality, we have
\[
\PP\left(\tS^{(n)}_k>\epsilon\right) 
\leq
\PP(\cB^{(n)}_k \text{ fails })
+ \epsilon^{-1} 
\EE
\left(\tS^{(n)}_k; \cB^{(n)}_k \right).
\]

By Lemma \ref{Blemma}, the first term tends to 0 as $k\to\infty$,
uniformly in $n$. For the second term, Lemmas \ref{differencelemma}
and \ref{EMlemma} combine to give
\begin{multline*}
\EE
\left(\tS^{(n)}_k; \cB^{(n)}_k \right)
\leq
cc_0(k+1)^{-\frac{1}{\alpha}+\frac{1}{2}+\delta}
+
cc_0\sum_{r=k+2}^{n^2} r^{-\frac{1}{\alpha}-\frac{1}{2}+\delta}
\\
+
cc_1(k+1)^{\frac{1}{2}}a_{n^2}^{-1}
+
cc_1 a_{n^2}^{-1} \sum_{c_2 n^2 \leq r \leq n^2} r^{-\frac{1}{2}}.
\end{multline*}
From the choice of $\delta$ and the fact that $a_{n^2}/n\to \infty$ as 
$n\to\infty$,
one obtains that 
all four terms on the RHS tend to 0 as $k\to\infty$ uniformly in $n$ 
such that $k\leq n^2$.

Hence indeed
$\PP\left(\tS^{(n)}_k>\epsilon\right)<\epsilon$ 
for all large enough $k$, uniformly in $n$
such that $k\leq n^2$, 
as required.$\hfill\Box$

\section{Path convergence}
\label{pathsection}

In this section we will state and prove results describing
the convergence of the distribution of the 
optimal paths for the discrete 
models to that for the limiting continuous model.

First we note that the optimal path for the continuous model
is well-defined:

\begin{proposition}
\label{uniquenessprop}
With probability 1,
there exists a unique $A^*\in\cC$ such that
$T=\sum_{i\in A^*} M_i$.
\end{proposition}

We will also define $A^{(n)*}$
as the set that achieves the maximum in
\[
T^{(n)}=\max_{A\in\cC^{(n)}}\sum_{i\in A} M^{(n)}_i.
\]
(or in the equivalent expression for $\tT^{(n)}$).
(Since the weight distribution is assumed to be continuous, 
this optimal set is almost surely unique).

It will be useful to extend the sequences
$(Y^{(n)}_i)_i$ and $(\tM^{(n)}_i)_i$ to all $i\in\NN$
(rather than only $i\leq n^2$); 
for example, we can put $\tM^{(n)}_i=0$, $Y^{(n)}_i=(0,0)$
for all $i\geq n^2$. We will consider always the
product topology when looking at convergence
of such infinite sequences.

We also use the product topology on $\cS$, the set of subsets of $\NN$.  
Thus, given $A$ and a sequence $A_k$ in $\cS$, we have $A_k\to A$ if, for
every $m$, $A_k \cap \{1,\dots,m\}$ is equal to $A \cap \{1,\dots,m\}$
for all large enough $k$.  One has easily that any sequence $A_k$ has
at least one limit point, and also that if $A_k\in\cC$ for each $k$
then every limit point is also in $\cC$ (since if the limit point
contains $i$ and $j$, then $i$ and $j$ are in $A_k$ for some $k$, and
hence $Y_i\sim Y_j$).

The following theorem is our first path convergence result. 
Later (in Theorem \ref{paththeorem2}) we will use 
it to prove a more direct convergence result
concerning the optimal paths viewed as random subsets of $[0,1]^2$.

\begin{theorem}
\label{paththeorem}
$\Big(\big(Y^{(n)}_i\big)_{i\in\NN}, A^{(n)*}\Big)
\to
\Big(\big(Y_i\big)_{i\in\NN}, A^{*}\Big)$
in distribution as $n\to\infty$.
\end{theorem}

Before proving Proposition \ref{uniquenessprop} and 
Theorem \ref{paththeorem}, we need the following fact:
\begin{lemma}
\label{pathlemma}
With probability 1, the following holds:
if $A_j$ is a sequence in $\cC$ converging to
a limit $A$, then
$\lim_{j\to\infty}\sum_{i\in A_j} M_i=\sum_{i\in A} M_i$.
\end{lemma}

\noindent
[N.B.\ this result is not true in general for sequences $A_j$ in $\cS$ 
(unless $\alpha<1$ so that $\sum M_i<\infty$ a.s.)]
\medskip

\noindent\textit{Proof:}
If $A_j \cap \{1,\dots,m\}=A \cap \{1,\dots,m\}$ then
\[ \left|\sum_{i\in A_j} M_i - \sum_{i\in A} M_i \right|
\leq \sup_{\tA\in\cC} \sum_{i\in\tA,i>m} M_i =S_m. \]
But $S_m\to 0$ a.s.\ as $m\to\infty$ (from Lemma \ref{Slemma}),
and $A_j \cap \{1,\dots,m\}=A \cap \{1,\dots,m\}$
eventually for all $m$, so we are done.$\hfill\Box$

\medskip

\noindent\textbf{Proof of Proposition \ref{uniquenessprop}}:
Recall that 
\[
T_k=\sup_{A\in\cC}\sum_{i\in A, i\leq k} M_i.
\]
Since the sum on the RHS depends only 
on the intersection of $A$ with $\{1,\dots,k\}$,
we only need to consider the max over finitely 
many $A\subseteq\{1,\dots,k\}$. 
Thus there exists some $A^*_{k}$ which achieves the sup.


We consider the sequence $A^*_k$. 
As observed above, this sequence has
at least one limit point $A^*\in\cC$, 
and by Lemma \ref{pathlemma},
$\sum_{i\in A^*} M_i=\lim_{k\to\infty} \sum_{i\in A^*_k} M_i
=\lim_{k\to\infty} T_k = T$.

Now we wish to show that in fact a unique $A^*$ 
achieves the sum $T$.

Suppose instead that there are two such optimising sets in $\cC$.
Then there is some $k$ that is contained in one but not the other.

In that case,
\begin{align*}
\sup_{A\in\cC, k\notin A} \sum_{i\in A} M_i
&=
\sup_{A\in\cC, k\in A} \sum_{i\in A} M_i
\\
&=M_k+\sup_{A\in\cC, k\in A} \sum_{i\in A, i\neq k} M_i,
\end{align*}
which gives
\begin{equation}\label{event}
M_k=
\sup_{A\in\cC, k\notin A} \sum_{i\in A} M_i
-
\sup_{A\in\cC, k\in A} \sum_{i\in A, i\neq k} M_i.
\end{equation}
We will show that this event has probability 0 for each $k$; 
then, by countable additivity, we are done.

The RHS of (\ref{event}) does not depend on $M_k$;
in fact, it is a function of the collection
$\big( (M_i)_{i\in\NN\setminus\{k\}}, (Y_i)_{i\in\NN} \big)$.
We condition on the value of this collection.
Then the RHS is a constant, while the random variable
$M_k$ on the LHS has a continuous distribution.
(Specifically, the distribution
of $M_k^{-\alpha}$ conditional on this collection
is uniform on the interval 
$\left( M_{k-1}^{-\alpha}, M_{k+1}^{-\alpha}\right)$,
since the sequence $\left(M_i^{-\alpha}\right)$ 
forms the points of a Poisson process).
Hence the event (\ref{event}) has probability 0, 
as required.$\hfill\Box$

\medskip

\noindent\textbf{Proof of Theorem \ref{paththeorem}}:
We will use the same couplings as in the proof of
Theorem \ref{Ttheorem}.
As at (\ref{conv}), we then have that 
$\PP\left(
\cC_{k_n}^{(n)}\neq \cC_{k_n}
\right)
\to 0$ as $n\to\infty$,
and that all of the quantities
$\sum_{i=1}^{k_n} \left| M_i - \tM^{(n)}_i \right|$,
$\sum_{i=1}^{k_n} \left\| Y_i - Y^{(n)}_i \right\|$,
$\tS^{(n)}_{k_n}$ and 
$S_{k_n}$ converge to 0
in probability as $n\to\infty$.

To prove Theorem \ref{paththeorem}, it will then suffice
to show in addition that for any $m$,
\[
\PP\left(A^{(n)*}\cap\{1,\dots,m\}\neq 
A^{*}\cap\{1,\dots,m\}\right)\to 0
\]
as $n\to\infty$.
For this, it's in turn enough to show that for all $r$,
\begin{gather}
\label{firstone}
\PP\left(r\in A^*,r\notin A^{(n)*}\right)\to 0
\\
\label{secondone}
\PP\left(r\notin A^*,r\in A^{(n)*}\right)\to 0
\end{gather}
as $n\to\infty$.
We will show (\ref{firstone}); an analogous argument gives 
(\ref{secondone}).

Define $T_{(-r)}=\sup_{A\in\cC, r\notin A} \sum_{i\in A} M_i$.

Suppose $r\in A^*$. 
Then $T_{(-r)} < T$ strictly.
(Otherwise, there is a 
sequence of members of $\cC$, none of which contain $r$,
whose weight converges to $T$;
then (by Lemma \ref{pathlemma})
this sequence has some limit point, 
itself a member of $\cC$ not containing $r$,
which attains the weight $T$. But this contradicts the 
uniqueness of $A^*$ established in Proposition
\ref{uniquenessprop}).

If $\cC^{(n)}_{k_n}=\cC_{k_n}$, then
\begin{align}
\nonumber
\max_{A\in\cC^{(n)}, r\notin A} \sum_{i\in A} \tM^{(n)}_i
&\leq
\max_{A\in\cC^{(n)}_{k_n}, r\notin A} \sum_{i\in A} \tM^{(n)}_i
+\tS^{(n)}_{k_n}
\\
\nonumber
&\leq
\max_{A\in\cC_{k_n}, r\notin A} \sum_{i\in A} \tM^{(n)}_i
+\sum_{i=1}^{k_n} \left[\tM^{(n)}_i - M_i\right]_+
+\tS^{(n)}_{k_n}
\\
\label{ineq1}
&\leq
T_{(-r)} + 
\sum_{i=1}^{k_n} \left[ \tM^{(n)}_i -M_i \right]_+
+ \tS^{(n)}_{k_n},
\\
\intertext{and similarly}
\label{ineq2}
\max_{A\in\cC^{(n)}, r\in A} \sum_{i\in A} \tM^{(n)}_i
&\geq
T -
\sum_{i=1}^{k_n} \left[ M_i - \tM^{(n)}_i  \right]_+
- S_{k_n}.
\end{align}
If also $r\notin A^{(n)*}$, then 
\[
\max_{A\in\cC^{(n)}, r\notin A} \sum_{i\in A} \tM^{(n)}_i
\geq
\max_{A\in\cC^{(n)}, r\in A} \sum_{i\in A} \tM^{(n)}_i,
\]
and using (\ref{ineq1}) and (\ref{ineq2}) we get 
\[
T-T_{(-r)}
\leq
\sum_{i=1}^{k_n} \left| M_i - \tM^{(n)}_i  \right|
+ S_{k_n} + \tS^{(n)}_{k_n}
\]
So altogether we obtain
\begin{multline*}
\PP\left(r\in A^*, r\notin A^{(n)*}\right)
\\
\leq
\PP\left(
\cC_k^{(n)}\neq \cC_k
\right)
+
\PP\left(
0<
T-T_{(-r)}
\leq
\sum_{i=1}^{k_n} \left| M_i - \tM^{(n)}_i  \right|
+ S_{k_n} + \tS^{(n)}_{k_n}
\right).
\end{multline*}
We have already observed above that the first probability on the RHS
tends to 0 as $n\to\infty$. The same is true for the second probability
on the RHS,
since all the terms on the right of the inequality converge 
to 0 in probability as $n\to\infty$, while the term in the middle of
the inequality does not depend on $n$.
Hence $\PP\left(r\in A^*, r\notin A^{(n)*}\right)\to 0$ as $n\to\infty$,
as required.\enpf

We now turn to the convergence of the paths
regarded as subsets of $[0,1]^2$.

First let $U^*=\bigcup_{i\in A^*} Y_i \cup \{(0,0), (1,1)\}$,
and take its closure $\Ubar^*$.

We expect that $\Ubar^*$ is connected with
probability 1, but we don't have a proof. To work around
this, we will use the fact that, at least,
there is a.s.\ a unique way to extend 
$\Ubar^*$ to a connected set while preserving the increasing
path property. (Here the increasing path property of a set 
means that if $y$ and $y'$ are two elements of the set then $y\sim y'$).

To see this, first note that if $\Ubar^*$ does ``contain jumps'',
then none of these jumps can span a rectangle of non-zero area.
That is, with probability 1 there is no rectangle $R$ 
of non-zero area such that $y\sim y'$ for all $y\in\Ubar^*$ 
and $y'\in R$.

For if there were, then $R$ would certainly contain some 
points $Y_j$; such $j$ could be added to $A^*$, 
increasing the weight of the path by $M_j$; this contradicts
the maximality of $A^*$.

So any jumps in $\Ubar^*$ consist only of horizontal or
vertical line segments. These segments can all be added to $\Ubar^*$
while still preserving the increasing path property, and
this gives a connected set. Conversely, any connected increasing
set containing $\Ubar^*$ must ``fill in'' these jumps.

Thus, define $P^*$ by setting $y\in P^*$ if:
\begin{itemize}
\item[(i)]$y\in\Ubar^*$, or
\item[(ii)]there exists $y',y''\in\Ubar^*$ with either
\begin{itemize}
\item[(a)]$y'(1)=y(1)=y''(1)$, $y'(2)<y(2)<y''(2)$, or
\item[(b)]$y'(2)=y(2)=y''(2)$, $y'(1)<y(1)<y''(1)$.
\end{itemize}
\end{itemize}

This set $P^*$ (which we conjecture to be equal to 
the closure of $\bigcup_{i\in A^*} Y_i$ w.\ p.\ 1)
provides the distributional limit we need.

For each $n$, we define the object representing
the optimal path in the discrete problem indexed by $n$ as follows: 
order the points $\{Y_i^{(n)}, i\in A^{(n)*}\}$
in increasing order and join successive points by a straight line
(horizontal or vertical, of length $1/n$). Call the 
resulting path $P^{(n)*}$.

$P^*$ and $P^{(n)*}$ are regarded as subsets of $[0,1]^2$
and we use the Hausdorff metric:
\[
d_H(P_1,P_2)=
\sup_{x\in P_1}
\inf_{y\in P_2}
|x-y|
+
\sup_{x\in P_2}
\inf_{y\in P_1}
|x-y|.
\]

\begin{theorem}
\label{paththeorem2}
$P^{(n)*}\to P^*$ in distribution as $n\to\infty$.
\end{theorem}

\noindent\textit{Proof:}
Choose a probability space 
on which the convergence in Theorem \ref{paththeorem}
occurs almost surely. We will show that in this case
$P^{(n)*}\to P^*$ a.s.\ also. 

First consider any point $y\in P^*$. We will show 
that for all sufficiently large $n$ there is a point 
of $P^{(n)*}$ within distance $\epsilon/2$ of $y$.

There are two cases to consider.

First, suppose that $y$ is a limit 
of some sequence of points $Y_i\in A^*$. Choose
some $Y_i$ which is with distance $\epsilon/4$.
For large enough $n$, we have $Y_i^{(n)}\in A^{(n)*}$
and $|Y_i^{(n)}-Y_i|<\epsilon/4$ and we are done.

Otherwise, $y$ is on a vertical or horizontal line
between two points that \textit{are} limits
of sequences $Y_i\in A^*$. Call these endpoints $y^-$ and $y^+$.
For large enough $n$, there are $i^-$ and $i^+$ such that
$Y_{i^-}, Y_{i^+}\in P^{(n)*}$ with 
$|Y_{i^-}^{(n)}-y^-|<\epsilon/2$ and 
$|Y_{i^+}^{(n)}-y^+|<\epsilon/2$, as above.
Then by the increasing path property,
the subpath of $P^{(n)*}$ joining $Y_{i^-}$ and $Y_{i^+}$
passes within $\epsilon/2$ of every point on the
line segment joining $y^-$ to $y^+$,
and hence in particular within $\epsilon/2$ of $y$ as required.

Now for some $m\in\NN$ set $\epsilon=1/m$
and consider an increasing sequence of points 
$(0,0)=y^{(0)}, y^{(1)}, \dots, y^{(2m)}=(1,1)\in P^*$
such that the $L_1$ distance $d_1(y^{(j)}, y^{(j+1)})$ between
successive points is exactly $\epsilon$ for all $j$. 
(This is possible since $P^*$ is an increasing path and connected).

Now for large enough $n$ there is a sequence 
$\ty^{(0)}, \ty^{(1)}, \dots, \ty^{(2m)}$
of points of $P^{(n)*}$ such that
$d_1(\ty^{(j)}, y^{(j)})<\epsilon/2$ for all $j$. 
Then necessarily $\ty^{(0)}, \ty^{(1)}, \dots, \ty^{(2m)}$
is itself an increasing sequence 
and $d_1(\ty^{(j)}, \ty^{(j+1)})<2\epsilon$ for all $j$.

Using the increasing path property for $P^*$, 
we have that every point of $P^*$ is within $L_1$
distance $\epsilon/2$ of one of the $y^{(j)}$. 
Then since each $y^{(j)}$ is within $L_1$ 
distance $\epsilon/2$ of a point of $P^{(n)*}$,
we have that every point of $P^*$ is within $\epsilon$
of $P^{(n)*}$. 

Similarly, the increasing path property for $P^{(n)*}$
gives that every point of $P^{(n)*}$ is within 
$L_1$ distance $\epsilon$ of one of the $\ty^{(j)}$,
and thus in turn within distance $3\epsilon/2$ of 
some point of $P^*$.

Hence $d_H(P^*, P^{(n)*})<5\epsilon/2$,
for all large enough $n$. This works for any $\epsilon=1/m$,
so $P^{(n)*}\to P^*$ as required.\enpf

\section{Stable process directed percolation}
\label{stablesection}

In this section we consider a directed last-passage
percolation model based on stable L\'evy processes.
This is the stable version of the Brownian
directed percolation problem considered in \cite{OcoYor, hmo}. 

For $n\in\NN$, $t>0$, consider the random variable
\[ L(n,t) = \sup_{0=t_0 \leq t_1 \leq...\leq t_n=t} \sum_{i=1}^n
S^i_{t_{i-1}t_i}, \]
where $S^i_{st} = S^i_t-S^i_s$, and $S^i$ are i.i.d.\ $\alpha$-stable
processes for some $\alpha\in(0,2)$.
The Brownian version of this problem, in which the
stable processes are replaced by Brownian motions, 
gives a representation
for the largest eigenvalue process in 
``Hermitian Brownian motion''
(a matrix-valued process whose marginal at any fixed time has the 
GUE distribution) 
and has been much studied in various contexts
(see for example \cite{GUEs, GraTraWid, hmo, OcoYor}). We do
not have a random matrix interpretation of this stable process
version; however, an interesting connection could be
to the case of Wigner random matrices with heavy-tailed
entries considered by Soshnikov in \cite{Soshnikovheavy},
where a scaling is obtained for the largest eigenvalues 
which corresponds to the one we have 
observed for the heavy-tailed last-passage percolation problem.

We will show that the asymptotic behaviour of the distribution
of $L(n,t)$, as $n$ becomes large, is again
described by our continuous heavy-tailed last-passage
directed percolation problem.
Note that by scaling $L(n,tn) = t^{1/\alpha} L(n,n)$ in
distribution and hence we can just consider $L(n,n)$.

The processes $S^i$ have jump measure 
$c_+ x^{-\alpha-1}I_{x>0} +c_-|x|^{-\alpha-1}I_{x<0}$,
for some $c_+>0$ and $c_-\geq 0$.
The jumps play the role of weights for the percolation problem.

\begin{theorem}
\[
\left(\frac{\alpha}{c_+}\right)^{1/\alpha}n^{-2/\alpha}L(n,n)\to T
\] 
in distribution as $n\to\infty$, where $T$ is the last-passage time
in the continuous last-passage percolation model with index $\alpha$,
defined at (\ref{Tdef}).
\end{theorem}

A short argument is available in the case $\alpha<1$
(making use of the fact that the sum of all positive weights 
is finite), and we give this first. 

\subsection{Case $\alpha<1$}
Let $M_1^{(n)}$, $M_2^{(n)}, \dots$
be the set of positive jumps of the processes 
$S^1$, $S^2,\dots,S^n$ on the interval $[0,n]$,
written in descending order.

From the form of the jump measure, we can regard the ordered sequence
of jumps as a Poisson random measure. Thus by a suitable
transformation we can write the sequence of jumps in terms of a Poisson process and 
have that for any $n$,
\begin{multline*}
\left(\frac{\alpha}{c_+}\right)^{1/\alpha}n^{-2/\alpha}
\Big(M_1^{(n)}, M_2^{(n)},\dots,M_k^{(n)},\dots\Big)
\\
\stackrel{d}{=}
\Big(W_1^{-1/\alpha},(W_1+W_2)^{-1/\alpha},
\dots, 
(W_1+\dots+W_k)^{-1/\alpha},
\dots\Big),
\end{multline*}
where $W_i$ are i.i.d.\ exponential random variables with mean 1.

Now let $L^{(n)+}_k$ be the maximal
weight of a path, if one ignores all the weights except the $k$
largest positive weights. Just as in the discrete case,
one can show that 
\[
\left(\frac{\alpha}{c_+}\right)^{1/\alpha}n^{-2/\alpha}L_k^{(n)+}
\to T_k
\]
in distribution as $n\to\infty$, 
and one also has $T_k\to T$ in distribution
as $k\to\infty$,
where $T_k$ and $T$ are the last passage times
for the continuous problem as defined before.

Now let $L^{(n)+}$ be the maximal weight
of a path, if one considers all the positive
weights but ignores all the negative ones. Then
\[
L^{(n)+} - L^{(n)+}_k
\leq
\sum_{r=k+1}^\infty M_r^{(n)}.
\]
Now the distribution of 
$n^{-2/\alpha}\sum_{r=k+1}^\infty M_r^{(n)}$
does not depend on $n$, and converges to $0$ 
in distribution as $k\to\infty$
(since the sum of all the positive weights is a.s.\
finite). So we have
\[
\left(\frac{\alpha}{c_+}\right)^{1/\alpha}n^{-2/\alpha}L^{(n)+}
\to T
\]
in distribution.

Now consider the optimal path attaining $L^{(n)+}$.  Consider the sum
of (the absolute values of) all the negative weights along the path;
call it $S^{(n)-}$.  Since the positive and negative weights occur
independently, $S^{(n)-}$ has just the same distribution as the sum of
the negative weights for a single stable process between times $0$ and
$n$.  This is finite and on the scale $n^{1/\alpha}$ (in fact, the
distribution of $n^{-1/\alpha}S^{(n)-}$ is independent of $n$).  So
certainly $n^{-2/\alpha}S^{(n)-}\to 0$ in distribution as
$n\to\infty$.  Since $L^{(n)+}-S^{(n)-}\leq L(n,n)\leq L^{(n)+}$, we
obtain
\[
\left(\frac{\alpha}{c_+}\right)^{1/\alpha}n^{-2/\alpha}L(n,n)
\to T
\]
in distribution as $n\to\infty$, as required.

\subsection{Case $1\leq \alpha<2$}

\noindent\textit{Lower bound:}

As before,
\[
\left(\frac{\alpha}{c_+}\right)^{1/\alpha}n^{-2/\alpha}L_k^{(n)+}
\to T_k
\]
in distribution as $n\to\infty$.

Now consider a path realising $L_k^{(n)+}$ in this way;
(for definiteness, say the first such path in the
lexicographic order).

Let $\tT^{(n)}_k$ be the total weight of this path,
including all weights, and let 
$\tS=\tT^{(n)}_k-L_k^{(n)+}$.

The distribution of $\tS$ can be described as follows. 
Generate the $n$ independent processes $S^1,\dots, S^n$
from time $0$ to time $n$. Remove 
the $k$ largest positive jumps that occur in 
the $n$ processes in time $[0,n]$. Then $\tS$ has
the distribution of the altered value of $S^1_n$,
after the $k$ largest jumps from the set of processes
have been removed.

However, as $n\to\infty$, the probability 
that any of the $k$ largest jumps occur in the process
$S_1$ tends to 0 (it is no larger than $k/n$); so with
high probability, the procedure in the previous paragraph
does not alter the value of $S^1_n$.
Thus the limit in distribution of $n^{-1/\alpha}\tS$
is the distribution of $n^{-1/\alpha}S^1_n$
(which is independent of $n$). In particular,
$n^{-2/\alpha}\tS\to 0$ in probability.
Thus, for any $k$,
\begin{align*}
\left(\frac{\alpha}{c_+}\right)^{1/\alpha}n^{-2/\alpha}\tT^{(n)}_k
&=
\left(\frac{\alpha}{c_+}\right)^{1/\alpha}n^{-2/\alpha}\tS+
\left(\frac{\alpha}{c_+}\right)^{1/\alpha}n^{-2/\alpha}L_k^{(n)+}
\\
&\to T_k
\end{align*}
in distribution, as $n\to\infty$.
Since $L(n,n)\geq \tT^{(n)}_k$ and $T_k\to T$ as $k\to\infty$,
this establishes that $T$ is a lower bound 
for the limit in distribution of 
$\left(\frac{\alpha}{c_+}\right)^{1/\alpha}n^{-2/\alpha}L(n,n)$.

\medskip

\noindent\textit{Upper bound:}

We first need a lemma on the tail behaviour
of the difference between the supremum and infimum
of the stable process on an interval:

\begin{lemma}\label{jumpuplemma}
Let $S$ be a stable process with index $\alpha$ 
and jump measure 
$c_+ x^{-\alpha-1}I_{x>0} +c_-|x|^{-\alpha-1}I_{x<0}$,
where $c_+>0$ and $c_-\geq 0$.
Then 
\[
\PP\Big(\sup_{0\leq s\leq t\leq 1}\{S_t-S_s\} > x\Big)
\sim \frac{c_+}{\alpha} x^{-\alpha}
\]
as $x\to\infty$. That is,
the quantity $\sup_{0\leq s\leq t\leq 1}\{S_t-S_s\}$
has the same positive tail behaviour
as the size of the largest positive jump of $S$
in $[0,1]$
(which is also the same as the upper tail of 
$S_1$ and of $\sup_{0\leq t\leq 1}S_t$).
\end{lemma}

\noindent\textit{Proof:}
Let the running infimum and supremum processes be denoted by $I_t =
\inf\{S_s:0\leq s\leq t\}$ and $S^*_t = \sup\{S_s: 0\leq s\leq t\}$
respectively. Consider the reflected process $X_t = S_t - I_t$. Our
aim is to determine the tail behaviour of $X^*_1 = \sup\{X_t:0\leq
t\leq 1\}$.  
By Bertoin \cite{Bertoinbook}~Prop~VI.3, 
we know that for fixed $t$, the distribution of $X_t$
is the same as that of $S^*_t$. By standard results
(e.g. \cite{Bertoinbook}~Prop~VIII.4) we have 
\begin{equation}
\PP(X_1>x) \sim \PP(S^*_1>x) \sim \frac{c_+}{\alpha} x^{-\alpha},
\label{eq:tailasymp}
\end{equation}
as $x\to\infty$.

Thus we just need to show that $X^*_1$ has the same tail as $X_1$. 
The proof is analogous to that of
\cite{Bertoinbook}~Prop.~VIII.4, 
and we reproduce the argument 
here. An easy consequence of (\ref{eq:tailasymp}) is that
\[ 
\liminf_{x\to\infty} \PP(X^*_1>x) x^\alpha \geq
\frac{c_+}{\alpha}. 
\]
Now fix $\epsilon>0$ and note that the reflected process is Markov by
\cite{Bertoinbook}~Prop~VI.1. As the stable process scales in that 
$S_{\lambda t} \stackrel{d}{=} \lambda^{1/\alpha} S_t$, 
this property will be inherited by the
infimum and hence the reflected process itself, giving 
$X_{\lambda t} \stackrel{d}{=} \lambda^{1/\alpha} X_t$.
Since $S_t-X_t=I_t$ is decreasing in $t$, we also 
have that $X_1-X_t\geq S_1-S_t$ for all $t<1$. 
Applying these properties and denoting
by $\tau_x$ the first hitting time of the interval $(x,\infty)$ by the
reflected process $X$, we have
\begin{align*}
\PP\big(X_1>(1-\epsilon)x\big)
&\geq
\PP\big(X_1^*>x, X_1>(1-\epsilon)x\big)
\\
&\geq
\int_0^1 \PP(\tau_x\in dt)\PP(X_1-X_t>-\epsilon x)
\\
&\geq
\int_0^1 \PP(\tau_x\in dt)\PP(S_1-S_t>-\epsilon x)
\\
&=
\int_0^1 \PP(\tau_x\in dt)\PP(S_{1-t}>-\epsilon x)
\\
&\geq
\int_0^1 \PP(\tau_x\in dt)\PP(S_1>-\epsilon x)
\\
&=\PP(X_1^*>x)\PP(S_1>-\epsilon x).
\end{align*}
%
%
As $\PP(S_1>-\epsilon x) \to 1$ as $x\to\infty$ we have that
\[ 
\limsup_{x\to\infty} P(X^*_1>x) x^{\alpha} \leq
(1-\epsilon)^{-\alpha} c_+/\alpha, 
\]
and, as $\epsilon$ is arbitrary, we have the result.\enpf

Now for $1\leq i, j \leq n$, define
\[
X(i,j)=\sup_{j-1\leq s\leq t\leq j} \{S^i_t-S^i_s\}.
\]
Let $T(n,n)$ be the maximal passage time for the discrete
model with weights $X(i,j)$.
Then one can see by a direct sample path comparison that
$L(n,n)\leq T(n,n)$.

Applying Lemma \ref{jumpuplemma}, 
$a_N\sim \left(\frac{\alpha}{c_+}\right)^{1/\alpha}N^{1/\alpha}$,
where $a_N=\inf\{x:\PP(X(i,j)>x\leq 1/N\}$.
Thus, applying Theorem \ref{Ttheorem} for the discrete model,
\[
\left(\frac{c_+}{\alpha}\right)^{1/\alpha}n^{-2/\alpha} T(n,n)
\to T
\]
in distribution as $n\to\infty$.
This gives the required upper bound in distribution 
for the limit of $L(n,n)$.

\section{The case $\alpha=0$: 
convergence to the greedy path}
\label{greedysection}

In this section we consider the 
discrete last-passage percolation model in the
case $\alpha=0$. The distribution
$F$ is said to have a 
\textit{slowly varying tail}: for all $t>0$, 
\[
\frac{1-F(tx)}{1-F(x)}\to 1 \textit{ as } x\to\infty;
\]
equivalently, 
for all $s<1$,
\begin{equation}
\label{first}
\frac{F^{-1}(1-sv)}{F^{-1}(1-v)}\to \infty \text{ as } v\downarrow0.
\end{equation}

Now it is no longer possible to find a non-degenerate 
limit in distribution for $T^{(n)}$ as we did 
in Theorem \ref{Ttheorem}. 
Let $M^{(n)}_1$ be the maximum of an i.i.d.\ sample from $F$
of size $n^2$ as before. Let $b_n$ be any sequence of constants.
Then any limit point in distribution of the sequence 
$b_n^{-1} M^{(n)}_1$ must be concentrated on the set $\{0,\infty\}$.
From Proposition \ref{greedyprop} below, the same is true
if we replace $M^{(n)}_1$ by $T^{(n)}$.

However, convergence in distribution of the optimal
paths can still be obtained. In fact, the form
of the limiting distribution has a particularly
simple description in terms of a ``greedy algorithm''.
We give a multifractal analysis of this limiting object
in Section \ref{greedysubsection}.

The limiting object is defined as follows.
Given the locations $Y_1, Y_2,\dots$ i.i.d.\ uniform on $[0,1]^2$,
let $\cC$, the set of increasing paths, be defined as at (\ref{cCdef}).
We now define the 
\textit{greedy path} $A^*=A^*(Y_1, Y_2, \dots)\in \cC$ 
recursively as follows.
Let $1\in A^*$ always, and then, 
given $A^*\cap \{1,\dots,r\}$,
let $r+1\in A^*$ if and only if $Y_{r+1}\sim Y_i$ for every
$i\in A^*$, $i\leq r$.
One can describe $A^*$ as the first member of $\cC$ in 
the lexicographic order.

The discrete problem is defined in terms of the locations
$(Y_i^{(n)})$ and weights $(M_i^{(n)})$ as before.
Write $A^{(n)*}$ for the optimal path as in Section
\ref{pathsection}. The following theorem
gives the convergence of these optimal paths to the greedy path:

\begin{theorem}
\label{slowtheorem}
$\Big(\big(Y^{(n)}_i\big)_{i\in\NN}, A^{(n)*}\Big)
\to
\Big(\big(Y_i\big)_{i\in\NN}, A^{*}\Big)$
in distribution as $n\to\infty$.
\end{theorem}

Exactly as in Section \ref{paththeorem},
one can also define $P^*$ and $P^{(n)*}$ 
to represent the optimal paths regarded as subsets of $[0,1]^2$,
and obtain the convergence in distribution of
$P^{(n)*}$ to $P^*$ as $n\to\infty$ (under the Hausdorff metric).
In fact, the situation is considerably simpler here;
one can simply define $P^*$ to be the closure of
$\bigcup_{i\in A^*} Y_i$, since by the results
of Section \ref{greedysubsection} this set is connected w.p.\ 1.

Theorem \ref{slowtheorem} will follow from the next proposition:
\begin{proposition}
\label{greedyprop}
For all $r$,
\[
\PP\left(M^{(n)}_r > \sum_{i=r+1}^{n^2} M^{(n)}_i \right)\to 1
\text{ as } n\to\infty.
\]
\end{proposition}

\noindent\textbf{Proof of Theorem \ref{slowtheorem}}:
Fix some $\epsilon>0$. As in Proposition \ref{coupledprop}, 
for $n$ large enough,
we can find a coupling such that, with probability at least $1-\epsilon$,
\begin{gather}
\nonumber
\sum_{i=1}^k \left\| Y_i - Y^{(n)}_i \right\| < \epsilon
\\
\intertext{and}
\label{thesame}
\cC_k^{(n)}=\cC_k
\end{gather}
Suppose that (\ref{thesame}) holds and also that, for all $r\leq k$,
$M^{(n)}_r > \sum_{i=r+1}^{n^2} M^{(n)}_i$.
Then indeed 
\begin{equation}
\label{thesame2}
A^{(n)*}\intersect\{1,\dots,k\} = A^*\intersect\{1,\dots,k\},
\end{equation}
where $A^*$ is the ``greedy path''. 
Then using Proposition \ref{greedyprop},
we can find $n$ such that
(\ref{thesame2}) holds 
with probability at least $1-2\epsilon$.
This gives the convergence in distribution in Theorem 
\ref{slowtheorem}.~$\hfill\Box$

\subsection{Proof of Proposition \ref{greedyprop}}

\begin{lemma}
\label{theme}
Let $\epsilon>0$ and $C>0$.
There exists $U(C,\epsilon)<1$ such that if 
$U(C,\epsilon)<u_1<(1+\epsilon)u_2-\epsilon$, then
\[
F^{-1}(u_1)<F^{-1}(u_2)(1-u_1)^{-C}(1-u_2)^{C}.
\]
\end{lemma}

\noindent\textit{Proof:}
Let $t=(1+\epsilon)^{2C}$.
From (\ref{first})
there exists $V>0$ such that for all 
$v<V$,
\[
\frac{F^{-1}\left(1-\frac{v}{1+\epsilon}\right)}{F^{-1}(1-v)}>t.
\]
Iterating,
\[
\frac{F^{-1}\left(1-\frac{v}{(1+\epsilon)^m}\right)}{F^{-1}(1-v)}>t^m,
\text{ for all }m\in\NN,
\]
and so in fact, for any $v_2<v_1<V$,
\[
\frac{F^{-1}(1-v_2)}{F^{-1}(1-v_1)}
>t^{\left\lfloor
\log_{1+\epsilon}\frac{v_1}{v_2}\right\rfloor}.
\]
Putting $u_1=1-v_1$, $u_2=1-v_2$ and $U(C,\epsilon)=1-V$,
we have that for $U(C,\epsilon)<u_1<u_2$,
\[
F^{-1}(u_2)>
t^{
\left\lfloor
\log_{1+\epsilon}\frac{1-u_1}{1-u_2}\right\rfloor}
F^{-1}(u_1).
\]
Now whenever $z\geq 1$, then $\lfloor z\rfloor\geq z/2$,
so restricting to $(1-u_1)>(1+\epsilon)(1-u_2)$
we get
\[
F^{-1}(u_2) > t^{
\frac{1}{2}\log_{1+\epsilon}\frac{1-u_1}{1-u_2}}
F^{-1}(u_1),
\]
which rearranges to the desired result.$\hfill\Box$

\begin{lemma}
\label{variation}
Fix $\epsilon>0$ and $C>0$. 
If $u_2\in(0,1)$ is sufficiently close to 1,
then for all $u_1$ with 
$0<u_1<(1+\epsilon)u_2-\epsilon$,
\[
F^{-1}(u_1)<F^{-1}(u_2)(1-u_1)^{-C}(1-u_2)^{C}.
\]
\end{lemma}

\noindent\textit{Proof:}
From Lemma \ref{theme}, we already know that this is true
when $u_1>U(C,\epsilon)$. Now if $u_1\leq U(C,\epsilon)$, 
then $F^{-1}(u_1)\leq F^{-1}(U(C,\epsilon))$. So it will suffice
to show that for $u_2$ sufficiently close to 1, the RHS is always 
at least $F^{-1}(U(C,\epsilon))$. In fact, we will show 
that the RHS tends to $\infty$ as $u_2\uparrow 1$, uniformly in $u_1$.

For any $u_2$, the RHS is minimised by $u_1=0$. 
So we wish to show that 
\begin{equation}
\label{needed}
F^{-1}(u_2)(1-u_2)^{C}
\to\infty \text{ as } u_2\uparrow1.
\end{equation}


Fix $\tu_1>U(2C,\epsilon)$. Then Lemma \ref{theme} gives,
for $u_2$ sufficiently close to 1,
\begin{align*}
F^{-1}(u_2)
&\geq F^{-1}(\tu_1)
(1-\tu_1)^{2C}(1-u_2)^{-2C}
\\
&=c(1-u_2)^{-2C}
\end{align*}
for some constant $c$.  This gives the desired convergence to $\infty$
in (\ref{needed}).~$\hfill\Box$

\begin{lemma}
\label{expectation}
Let the r.v.\ $U$ be uniform on $(0,1)$. 
\[ \lim_{u\to1}\frac{\EE(F^{-1}(U)\given U\leq u)}{F^{-1}(u)(1-u)}
=0. \]
\end{lemma}

\noindent\textit{Proof:}
Take any $C>0$ and $\epsilon>0$.
If $u$ is close enough to 1, then using Lemma \ref{variation},
\begin{align*}
\EE(F^{-1}(U); U\leq u) &= \int_0^u F^{-1}(v)dv \\
&\leq \int_0^{(1+\epsilon)u-\epsilon} F^{-1}(u)(1-v)^{-C}(1-u)^{C}dv +
\int_{(1+\epsilon)u-\epsilon}^u F^{-1}(u)dv \\
&= F^{-1}(u)\left\{(1-u)^{C}\frac{(1+\epsilon)(1-u)^{-C+1}}{C-1} +
\epsilon(1-u) \right\} \\
&=F^{-1}(u)(1-u)\left\{\frac{(1+\epsilon)^{-C+1}}{C-1}+\epsilon \right\}.
\end{align*}

We can now choose $\epsilon$ as small as desired, 
and then $C$ as large as desired, to give an upper 
bound on 
\[ \limsup_{u\to1}\frac{\EE(F^{-1}(U); U\leq u)}{F^{-1}(u)(1-u)} \]
which is arbitrarily close to 0.

Finally, for $u\geq 1/2$, 
\[ \EE(F^{-1}(U)\given U\leq u) \leq 2\EE(F^{-1}(U);U\leq u) \]
and the result follows.$\hfill\Box$

\medskip

\noindent\textbf{Proof of Proposition \ref{greedyprop}}:
We use the representation  
\[
\left(
M^{(n)}_1,\dots,M^{(n)}_{n^2}
\right)
=
\left(
F^{-1}(U^{(n)}_1,\dots, U^{(n)}_{n^2})
\right),
\]
where $(U^{(n)}_1,\dots, U^{(n)}_{n^2})$
are the order statistics of an i.i.d.\ sample of size $n^2$ from 
the uniform distribution on $(0,1)$, written
in decreasing order.

We need to show that
\[ \PP\left( \sum_{m=r+1}^{n^2} F^{-1}(U^{(n)}_m)\geq
F^{-1}(U^{(n)}_r)\right)\to 0 \text{ as } n\to\infty. \]
Let $R>0$ and suppose that $u>1-R/n^2$.
Then
\begin{align}
\nonumber
\frac{n^2\EE\big(F^{-1}(U)\given U\leq u\big)}{F^{-1}(u)}
&=\frac{R \EE\big(F^{-1}(U)\given U\leq u\big)}{F^{-1}(u)R/n^2} \\
\nonumber
&\leq \frac{R\EE\big(F^{-1}(U)\given U\leq u\big)}{F^{-1}(u)(1-u)} \\
\label{bound}
&\leq R\sup_{u>1-R/n^2}\frac{\EE\big(F^{-1}(U)\given U\leq u\big)}
{F^{-1}(u)(1-u)}.
\end{align}
Now for (almost) all $u$,
\[
\EE\left(\sum_{m=r+1}^{n^2} F^{-1}(U^{(n)}_m)
\bigg\vert
U^{(n)}_r=u\right)
=(n^2-r)\EE\left(F^{-1}(U)\given U\leq u\right),
\]
so for (almost) all $u>1-R/n^2$ we have, from Markov's inequality,
\begin{multline*}
\PP\left(\sum_{m=r+1}^{n^2} F^{-1}(U^{(n)}_m)\geq F^{-1}(M^{(n)}_r)
\bigg\vert U^{(n)}_r=u
\right)
\\
\begin{aligned}
&\leq
\frac{1}{F^{-1}(u)}
\EE\left(\sum_{m=r+1}^{n^2} F^{-1}(U^{(n)}_m)
\bigg\vert
U^{(n)}_r=u\right)
\\
&\leq
n^2
\frac{
\EE\left(F^{-1}(U)\given U\leq u\right)
}
{F^{-1}(u)}
\\
&\leq
R\sup_{u>1-R/n^2}
\frac{
\EE\big(F^{-1}(U)\given U\leq u\big)
}
{F^{-1}(u)(1-u)}
\end{aligned}\end{multline*}
using (\ref{bound}).

So in fact, integrating over $u>1-R/n^2$,
\begin{multline*}
\PP\left(\sum_{m=r+1}^{n^2} F^{-1}(U^{(n)}_m)\geq F^{-1}(M^{(n)}_r)
\bigg\vert U^{(n)}_r>1-R/n^2
\right)
\\
\leq
R\sup_{u>1-R/n^2}
\frac{
\EE\big(F^{-1}(U)\given U\leq u\big)
}
{F^{-1}(u)(1-u)}.
\end{multline*}
Then
\begin{multline*}
\limsup_{n\to\infty}
\PP\left(\sum_{m=r+1}^{n^2} F^{-1}(U^{(n)}_m)\geq F^{-1}(U^{(n)}_r)\right)
\\
\leq
\limsup_{n\to\infty}
\PP\left(U^{(n)}_r\leq 1-R/n^2\right)
+\limsup_{n\to\infty}
R\sup_{u>1-R/n^2}
\frac{
\EE\big(F^{-1}(U)\given U\leq u\big)
}
{F^{-1}(u)(1-u)}.
\end{multline*}
The second term is 0 by Lemma \ref{expectation}, for every $R$.
The first term is the probability $\PP(B_{n^2,R/n^2}\leq r-1)$,
where $B_{n^2,R/n^2}$ is a binomial $(n^2,R/n^2)$ random variable.
This converges as $n\to\infty$ to the probability
$\PP(P_R\leq r-1)$ where $P_R$ is a Poisson mean $R$ random 
variable. Thus this probability can be made as small as desired by choosing $R$
large, and hence the limsup is in fact $0$ as required.$\hfill\Box$

\subsection{The properties of the greedy path}
\label{greedysubsection}

\begin{figure}[t]
\centering
\epsfig{figure=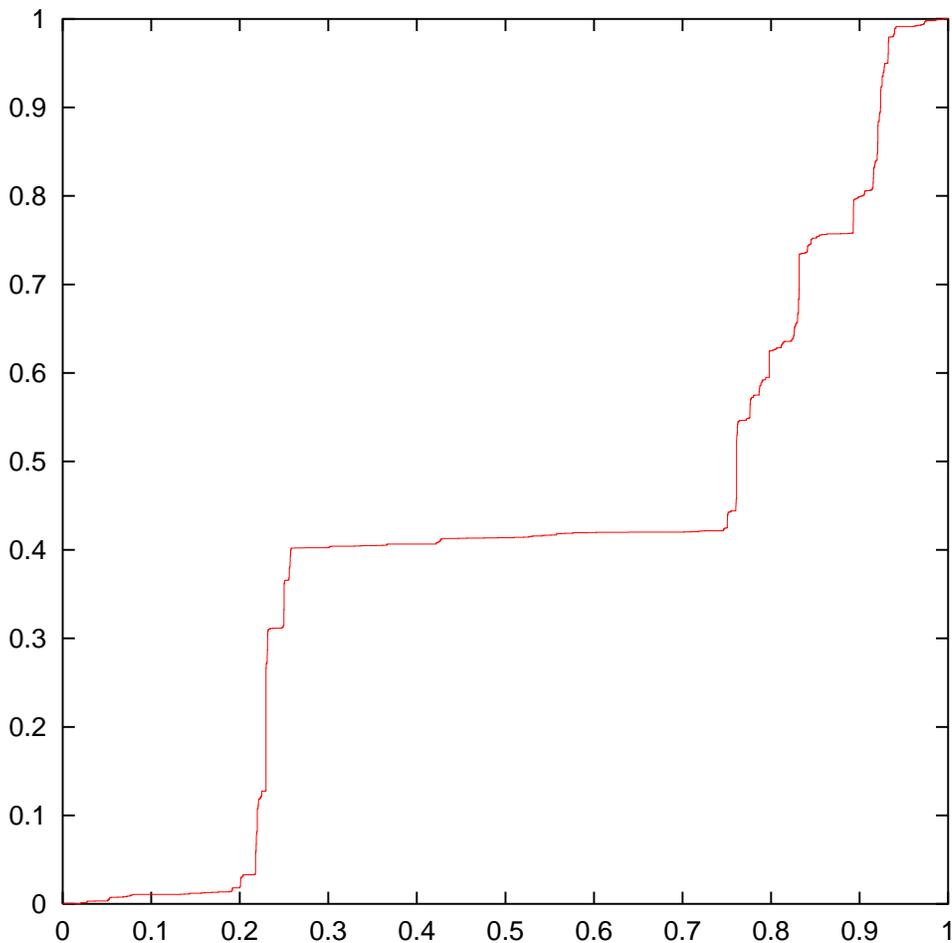, angle=270, width=1.3\linewidth}
\caption{A simulation from the distribution of the  ``greedy path'' 
which occurs as the limit of the distribution of optimal paths
in the case $\alpha=0$.
\label{greedyfig}}
\end{figure}

In this section we discuss the properties of the greedy path
which we have obtained as a distributional limit 
of the optimal path for the discrete problem when $\alpha=0$. 
The path can be regarded as a function $y=G(x)$ from $[0,1]$
to itself in a natural way 
(specifically, one could define
$G(x)=\sup\{y: Y_i=(x',y) \text{ for some } x'\leq x 
\text{ and some } i\in A^*\}$. 
See Figure \ref{greedyfig} for a realisation of the greedy path). 
This function $G$ is monotone non-decreasing, and hence defines
a measure $\mu$ on $[0,1]$. 
We will show that $\mu$ 
is a random self-similar measure which is singular with
respect to Lebesgue measure, and we will be able to compute its
multifractal spectrum.

We recall the definition of the multifractal spectrum for our
setting. Let $B_r(x)$
denote the ball of radius $r$ around the point $x\in \br$. 
For $a\geq 0$,
we define $E_a$, the set of
points at which the measure has local dimension $a$ by
\[ E_{a} = \left\{ x: \lim_{r\downarrow 0} \frac{\log
  \mu(B_r(x))}{\log r} = a \right\}. \]
The multifractal spectrum is then defined to be
\[ f(a) = \dim_H(E_{a}), \]
where $\dim_H(A)$ denotes the Hausdorff dimension of a set $A$. In our
setting we have 
\[ E_{a} = \left\{x\in [0,1]: \lim_{r\downarrow 0}
\frac{\log(G(x+r)-G(x-r))}{\log r} = a \right\}. \]

In particular we note that if $\mu$ had a density with respect to
Lebesgue measure, then the spectrum would be the function 
$f(a)=0$ for all $a\neq 1$ and $f(1)=1$.

There are a number of papers making rigorous the 
\textit{multifractal formalism}, the heuristic argument
for computing the multifractal spectrum in terms of the Legendre
transform of the moment measures, and we will be able to set our
measure in a framework within which we can apply this formalism. The
study of the multifractal spectrum for random self-similar measures is
the topic of \cite{fal, ap, bhj} where the underlying assumptions are
successively weakened.  

A \textit{scaling law} on a space $E$ consists of a probability space $(\Omega,
\cf, \bp)$ and for each $\omega\in \Omega$ a collection of weights and
maps $(\phi_1(\omega),p_1(\omega),\dots,\phi_N(\omega),p_N(\omega))$,
where $p_i\in \br_+$ and $\phi_i:E \to E$ is a contraction with
Lipschitz constant $r_i$. For a given scaling law a random
self-similar measure is 
a measure $\mu$ which satisfies the distributional equality
\[ \mu(\cdot) = \sum_{i=1}^N p_i \mu_i (\phi_i^{-1}(\cdot)), \]
where $\mu_i$ are i.i.d.\ copies of $\mu$ 
(independent of the weights and maps).
The support of the measure is typically a random self-similar set.

The multifractal formalism enables the
multifractal spectrum for the random self-similar measure to be
calculated in the following way. Let $m(q,\theta) = E( \sum_i p_i^q
r_i^{\theta})$ and $\beta(q) = \inf\{\theta:m(q,\theta)\leq 1\}$. 
Under the formalism the multifractal spectrum is the Legendre transform of
$\beta(q)$,
\begin{equation}
 f(a) = \inf_{q\in\br}\{a q + \beta(q)\}. \label{eq:multform}
\end{equation}
We will now give a more formal version.

We introduce a little notation. 
We write $\cT_n = \{1,\dots,N\}^n$ for
the sequences which index the sets after $n$ applications of the scaling
law, and $\cT$ for the tree 
$\cup_{n=0}^{\infty} \cT_n$. 
Let $(\Omega^{\otimes \cT},
\cf^{\otimes \cT}, \bp^{\otimes\cT})$ denote the product probability space for
random variables on the tree; for each node,
we have an independent copy of the scaling law.
Now for each $\bfi\in\cT_n$, define
\begin{align*}
p_{\bfi} = 
p_{i_1}(\omega_{\emptyset})p_{i_2}(\omega_{i_1})\dots 
p_{i_n}(\omega_{i_1 i_2\dots i_{n-1}}) 
\\
\intertext{and}
r_{\bfi} = 
r_{i_1}(\omega_{\emptyset})r_{i_2}(\omega_{i_1})\dots 
r_{i_n}(\omega_{i_1 i_2\dots i_{n-1}}).
\end{align*}
The total mass of the random measure
over the unit interval is given by 
$W = \lim_{n\toi} \sum_{\bfi\in\cT_{n}} p_{\bfi}$. 
In our setting we will consider random
probability measures, so that $W=1$. 
The other limit random variable we need is
is $W(q) =\lim_{n\to\infty} W_n(q)$ where 
$W_n(q) = \sum_{\bfi\in\cT_n} p_{\bfi}^q r_{\bfi}^{\beta}$. 

We also define the set $I_\beta\subseteq\RR$ as follows:
$q*\in I_\beta$ if, for some $a\geq0$, 
the infimum $\inf_{q\in\RR}\{aq+\beta(q)\}$ is non-negative 
and is achieved at $q*$.

Finally the \textit{strong open set condition}
is that there is an open set $O$ such that,
with probability 1, one has that 
$\phi_i(O)\cap\phi_j(O)=\emptyset$ for $i\neq j$, that $\cup_{i=1}^N
\phi_i(O) \subset O$ and that $\mu(O)>0$. 

We now state a
version of the main result of \cite{bhj} which can be applied in our
setting. 

\begin{lemma}\label{lem:mf}
Let $\mu$ be a random self-similar probability measure satisfying the 
strong open set condition. If the following three sets of
conditions are satisfied: 
\begin{enumerate}
\item $-\infty< \be \sum_i p_i\log p_i <0$;
\item For all $q\in I_\beta$, 
$\be \sum_i (\log p_i) p_i^q r_i^{\beta(q)}<\infty$ 
and 
$\be \sum_i (\log r_i) p_i^q r_i^{\beta(q)}<\infty$;
\item  
For all $q\in I_\beta$,
$\be 
\sum_i \left((\log p_i)^2 + (\log r_i)^2\right) p_i^q r_i^{\beta(q)}<\infty$
and 
$\be W(q) \log_+W(q)<\infty$; 
\end{enumerate}
then the multifractal formalism holds in the following sense.
Define 
\[ \beta^*(a) = \inf_{q\in \br}\{ a q + \beta(q)\}. \]
Then for any $a\geq0$,
\[ f(a) = \max \{ \beta^*(a), 0\},\]
with probability 1, and 
$E_a = \emptyset$ with probability 1 if 
$\beta^*(a)<0$. 
\end{lemma}

We now return to the measure arising from our greedy path. 
From the definition, we can describe the greedy path 
recursively as follows: choose a point $Y=(Y(1), Y(2))$ uniformly in 
the box $[0,1]^2$. Then the original path
is the union of two independent greedy paths, 
one scaled to lie in $[0,Y(1)]\times[0,Y(2)]$ 
and the other to lie in $[Y(1),1]\times[Y(2),1]$.
From this representation we can regard the induced measure
as a random self-similar measure for a scaling law.
Let $V$ and $\tV$ be independent uniform random variables.
Then the scaling law has 
$N=2$ with the set of weights
$(p_1,p_2)=(\tV,1-\tV)$ and the set of contractions $(\phi_1,\phi_2)$,
where $\phi_1(x) = Vx$ has contraction factor $r_1=V$ and $\phi_2(x) =
1-(1-V)x$ has contraction factor $r_2=1-V$.
Then the random self-similar measure $\mu$ satisfies
$\mu(\cdot) = \sum_{i=1}^2 p_i \mu_i(\phi^{-1}_i(\cdot))$ in distribution. 
From the construction it is clear that the measure is a probability
measure whose support is the unit interval. 

We note that our measure does not fit into the framework of \cite{fal}
or \cite{ap}. The interval does not satisfy the strong separation
condition as required in \cite{fal} (for a definition see \cite{fal})
but instead it satisfies the strong open set condition with
$O=(0,1)$. This ensures that the overlap between a pair of
contractions applied to the unit interval occurs at one point.
It also uses uniform random variables for the contraction ratios and
therefore there is no strictly positive lower bound on the contraction
ratios as required in \cite{ap}.

\begin{theorem}\label{thm:mfs}
The measure corresponding to the greedy path has for, any given $a\geq 0$,
with probability 1,
\[ f(a) = \left\{ \begin{array}{cl} \sqrt{8a} - a -1, &
  3-2\sqrt{2} \leq a \leq 3+2\sqrt{2}, \\ 0, & \mbox{ otherwise.}
  \end{array} \right. \] 
If $a\in[0,3-2\sqrt{2})\cup(3+2\sqrt{2},\infty)$ then $E_a=\emptyset$ 
with probability 1.
\end{theorem}

\noindent
{\it Proof:}
We determine the multifractal spectrum using the multifractal
formalism. We need to consider $m(q,\theta) =
E(\tV^q V^\theta + (1-\tV)^q (1-V)^\theta)$ and $\beta(q) =
\inf\{\theta:m(q,\theta)\leq 1\}$. It is straightforward to compute
these quantities and we have
\[ m(q,\theta) = 2\int_0^1 y^q dy \int_0^1 x^\theta dx =
\frac2{(1+q)(1+\theta)} \mbox{ for $q>-1, \theta>-1$}. \]
Thus
\[ \beta(q) = \frac2{q+1} -1, \]
and we can calculate that 
\[
\beta^*(a)=\sqrt{8a}-a-1
\]
for all $a\geq 0$.
We have that $\beta^*(a)$ is non-negative on the interval
$[3-2\sqrt2, 3+2\sqrt2]$ and negative elsewhere. 
Hence, if we can establish the three conditions of
Lemma~\ref{lem:mf}, to justify the formalism, we will have proved our
theorem.

Another calculation gives that $I_\beta=[1-\sqrt2,1+\sqrt2]\subset(-1/2,3)$,
and that for all $q\in I_\beta$, one also has $\beta(q)\in(-1/2,3)$.

For condition (1), we can compute $\be \sum_i p_i\log p_i = -1/2$.

For (2) straightforward calculations give  $\be \sum_i (\log p_i) p_i^q
r_i^{\beta} = 2((q+1)^2(1+\beta))^{-1}<\infty$ and $\be \sum_i (\log r_i) p_i^q
r_i^{\beta} = 2((q+1)(1+\beta)^2)^{-1}<\infty$, for all $q>-1, \beta>-1$.

Finally for the conditions (3) we have to do some work.
It is easy to calculate the first condition
\[ \be \sum_i \left((\log p_i)^2 + (\log r_i)^2\right) p_i^q
r_i^{\beta} = \frac{2}{(1+q)^3(1+\beta)} +
\frac{2}{(1+q)(1+\beta)^3}<\infty, \]
for $q>-1,\beta>-1$.

Thus we only have to verify the final condition.

We begin by observing that all we need is to prove that for some
$\epsilon>0$, $\be W(q)^{1+\epsilon}<\infty$ for each $q$. To do this
we will show that $\be W_n(q)^{1+\epsilon}$ converges.  First we
observe that $W_n(q)$ is a martingale by the definition of $\beta(q)$ 
and we compute the bracket process
\begin{eqnarray*} 
[W(q)]_n &=& 
\sum_{i=1}^n \be \Big(\big(W_i(q)-W_{i-1}(q)\big)^2 | \cf_{i-1}\Big) 
\\
 &=& \sum_{i=1}^n \be \left(\left(\sum_{j\in \cT_{i-1}} p_j^qr_j^{\beta(q)}
 (\chi_j-1)\right)^2\Big| \cf_{i-1}\right), 
\end{eqnarray*}
where $\chi_j = \tV_j^qV_j^{\beta(q)} + (1-\tV_j)^q(1-V_j)^{\beta(q)}$ and the
$(V_j,\tV_j)$ are independent of each other and independent over $j$. By
definition of $\beta$ we know that $\be \chi_j = 1$ and we can also
compute 
\[ \be \chi_j^2 = \frac{2}{(1+2q)(1+2\beta)} +
2\frac{\Gamma(q+1)^2}{\Gamma(2q+2)}
\frac{\Gamma(\beta+1)^2}{\Gamma(2\beta+2)}. \]
One can easily check that this quantity is finite over 
$q$ and $\beta$ in $(-1/2,3)$, and hence for all $q\in I_\beta$.
Thus we have
\[ [W(q)]_n = \sum_{i=1}^n W_{i-1}(2q,2\beta) (\be \chi^2-1), \]
where we write $W_n(a,b) = \sum_{\bfi\in\cT_n} p_{\bfi}^a r_{\bfi}^b$.

With the bracket process we can control the moments of the martingale
as for any $\gamma\geq 1$ there is a constant $c_{\gamma}$ such that
\[ \be W_n(q)^\gamma \leq c_{\gamma} \be [W(q)]_n^{\gamma/2}. \]
Thus we need to compute the moments of the bracket process. As we will take
$1<\gamma<2$ and all terms in the sum are positive, straightforward estimates give  
\begin{eqnarray*}
\be [W(q)]_n^{\gamma/2} &=& \be(\sum_{i=1}^n W_{i-1}(2q,2\beta) (\be
\chi^2-1))^{\gamma/2} \\
&\leq& \be \sum_{i=1}^n W_{i-1}(2q,2\beta)^{\gamma/2} (\be
\chi^2-1)^{\gamma/2} \\
&\leq& \sum_{i=1}^n \be W_{i-1}(\gamma q,\gamma\beta) (\be
\chi^2-1)^{\gamma/2}.
\end{eqnarray*}
Thus we just need to find $\be W_i(\gamma q, \gamma\beta) = (\be(
\tV^{\gamma q}V^{\gamma\beta} + (1-\tV)^{\gamma
  q}(1-V)^{\gamma\beta}))^i$. Integration gives
\[ \be W_i(\gamma q, \gamma\beta) = \left(\frac{2}{(1+\gamma
  q)(1+\gamma\beta)}\right)^i = \left(\frac{2(1+q)}{(1+\gamma
  q)(1+q+\gamma(1-q))}\right)^i. \]
Thus our result will hold if we can establish that over the range of
  $q$, we can find a $\gamma>1$ such that
\[ \frac{2(1+q)}{(1+\gamma q)(1+q+\gamma(1-q))} < 1. \]
A numerical calculation with the quadratic formula shows that this is
the case and hence we have our result. \enpf

The explicit form of the spectrum shows that $f(a)>0$ for
$3-2\sqrt{2}<a<3+2\sqrt{2}$. We also observe that the set of
points for which $a=2$ has full dimension 1. The next result
shows that the measure corresponding to the path is ``singular'':

\begin{corollary}
With probability 1, the greedy path is continuous and strictly increasing, 
with zero derivative almost everywhere. 
\end{corollary}

\noindent
\emph{Proof:} 
These properties can be proved fairly directly from 
the construction of the greedy path, but here we 
deduce them immediately from the multifractal spectrum. 

By construction, the distribution of the path is symmetric in
the $x$ and $y$ coordinates; hence the continuity property 
and the property that the path is strictly increasing are equivalent. 

Any point of discontinuity of the path belongs to the set $E_0$,
by definition. But from Theorem \ref{thm:mfs} we have that
$E_0$ is a.s.\ empty, so the path is a.s.\ continuous as desired.

For the derivative, note that 
as $G$ is a distribution function
it is almost everywhere differentiable with non-negative derivative. 
Let $D_1$ denote the set of points where the path is differentiable
and has strictly positive derivative.
It is straightforward to see that $D_1\subset E_1$. 
Thus as $f(1) = \sqrt{8}-2<1$, $\bp-a.s.$, 
we have that the Lebesgue measure of $D_1$ is 0 with probability 1. \enpf

\section{Last-passage random fields and an Airy process}
\label{Airysection}


In this section we consider the extension of the results to the case
of a random field. To do this we give a Poisson random measure
construction of the continuous limit model.

Let $\mu$ denote a Poisson random measure on $\br^3_+$ with intensity
measure 
$\lambda(dx\,dy\,dz) = dx\,dy\,\alpha z^{-\alpha-1} dz$. That is, if
$N(x,y,z) = \mu((0,x)\times(0,y)\times (z,\infty))= \int_0^x \int_0^y
\int_z^{\infty} \mu(dx\,dy\,dz)$ denotes the
number of points in $(0,x)\times(0,y)\times (z,\infty)$, then this has
a Poisson distribution with mean $\be N(x,y,z) = xy z^{-\alpha}$,
and the number of points in disjoint sets are independent.

We can now extend our limiting model to this setting. To relate this
to our original model we can consider the unit square in $\br_+^2$ and
order the points of the Poisson random measure in decreasing order of their
$z$-coordinates and we recover the sequence of weights
$Z_i=M_i$ in the original model. 

We will write $(Y_i,Z_i)$ for a point of the Poisson random
measure, where the points are labelled as above in that we regard
$Y_i$ as the location in $\br^2_+$ and $Z_i$ as the weight in our continuous
last passage percolation model.  Let
\[ \cC_{xy} = \{A \subset \bn 
\mbox{ such that for all $Y_i \sim Y_j$ 
for all $i,j\in A$ and 
$Y_i \in [0,x]\times[0,y]$ for all $i\in A$}\} \]
and define
\[ T(x,y) = \sup_{A\in \cC_{xy}} \sum_{i\in A} Z_i. \]
It is clear by Theorem~\ref{Ttheorem} that this random variable will
exist for each fixed $x,y$. This can be extended 
to show that the random field
$\{T(x,y),x>0,y>0\}$ exists and arises as the limit of the last
passage model. Let 
\[ T^{(n)}(x,y) = \max_{\pi\in \Pi^n_{xy}} \sum_{v\in \pi} X(v), \]
where $\Pi^n_{xy}$ is the set of directed
paths from $(1,1)$ to $(\lceil nx \rceil, \lceil ny \rceil)$,
and let 
\[ \tilde{T}^{(n)}(x,y) = \frac{T^{(n)}(x,y)}{a_{n^2}}. \]

\begin{theorem}
The field $\{T(x,y), x>0,y>0\}$ 
exists almost surely and $\{\tilde{T}^{(n)}(x,y),
x>0,y>0\} \to \{ T(x,y), x>0,y>0\}$ in the sense of convergence of
finite dimensional distributions.
\end{theorem}

\noindent
\emph{Proof:}
Given a realisation of the PRM, the random variable
$T(x,y)$ is well defined for all $x,y$, and over any finite box
$[0,c_x]\times [0,c_y]$ we have $T(x,y)\leq T(c_x,c_y)$ whenever
$x\leq c_x, y\leq c_y$; for any fixed $c_x,c_y$ we have
$P(T(c_x,c_y)<\infty)=1$ and hence by countable additivity
\[ P( T(x,y)<\infty \text{ for all } x,y ) =1. \]



We now need to establish the convergence of finite dimensional distributions. 
For one particular point we already have the one-dimensional
convergence result given in Theorem~\ref{Ttheorem}. 
Exactly the same couplings 
between the discrete and continuous problems that we  
used to prove Theorem~\ref{Ttheorem} 
can be applied to extend the result to a finite collection of points
within any finite box. The size of the box is
arbitrary and we obtain convergence of all the finite dimensional
distributions. \enpf

We now proceed to define a stationary field on the whole of
$\br^2$. To do this we observe that by simple scaling $N(\lambda x, \nu y,
(\lambda\nu)^{1/\alpha} z) = N(x,y,z)$ in distribution. In particular 
we have
\begin{equation}
 T(\lambda x, \nu y) = (\lambda\nu)^{1/\alpha} T(x,y) \mbox{ in
  distribution.}\label{eq:scale}
\end{equation}
Now put 
\[ \Theta(u_1,u_2) = \exp\left(-\frac{(u_1+u_2)}{\alpha}\right) 
T(e^{u_1},e^{u_2})
\]
for all $u_1,u_2\in\RR$.
Then we have that for any $v_1,v_2\in\RR$,
the collections 
$\{\Theta(u_1,u_2), u_1,u_2\in\RR\}$
and 
$\{\Theta(u_1+v_1,u_2+v_2), u_1,u_2\in\RR\}$
have the same distribution.

The next two results concern the moments and correlations of 
this stationary field. 

\begin{proposition}\label{prop:fieldmoment}
For all $\beta\in(0,\alpha)$ and all $\bfu\in\br^2$, we have 
$\be \Theta(\bfu)^{\beta} = \be T(1,1)^{\beta}<\infty$.
\end{proposition}

\noindent
\emph{Proof:} This is a continuation  of the argument given to
establish the existence of $T(1,1)$ in the proof of
Lemma~\ref{Slemma}. Recalling the setting of that proof,
we have $T(1,1)=S_0\leq U_0$,
where for $k\geq 0$, we defined 
$U_k=\sum_{i=k+1}^\infty L_i (M_i-M_{i+1})$.
Here $L_i$ is the largest number of the first $i$ locations 
that can be included in an increasing path; 
there is a constant $c$ such that
$\EE L_i\leq c\sqrt i$ and $\EE L_i^2\leq ci$ for all $i$,
and the collections $(M_i)$ and $(L_i)$ are independent.

Since $T(1,1)\leq U_0$, it will be enough to show that $\EE U_0^\beta<\infty$.
If a finite collection of random variables each have a 
finite $\beta$th moment, then so does their sum.
Hence, since $U_0=L_1(M_1-M_2)+\dots+L_k(M_k-M_{k+1}) +U_k$,
it's enough to show both of the following:
\begin{itemize}
\item[(i)]for any $i$, $\EE\big[L_i(M_i-M_{i+1})\big]^\beta<\infty$;
\item[(ii)]for some $k$, $\EE U_k^\beta<\infty$.
\end{itemize}
For property (i), note that $\big[L_i(M_i-M_{i+1})\big]^\beta<(iM_1)^\beta$,
so it's enough to show that $\EE M_1^\beta<\infty$.
But $M_1=W_1^{-1/\alpha}$, where $W_1$ has exponential distribution
with mean 1. Thus $\EE M_1^\beta=\int_0^\infty w^{-\beta/\alpha}e^{-w}dw$,
which is finite for all $\beta<\alpha$ as required.

So it remains to show (ii).
Since $\beta<\alpha<2$, it is enough to show that $\EE U_k^2<\infty$
for all large enough $k$, and this is what we will do.

Recall that we can write $M_i=(W_1+\dots+W_i)^{-1/\alpha}$
where $W_j$ are i.i.d. exponential random variables
with mean 1. We also write $V_i=W_1+\dots+W_i$, so that $M_i=V_i^{-1/\alpha}$.
Then, using the fact that $(1+x)^{-1/\alpha}>1-x/\alpha$ for all $x>0$,
we have that for all $i$,
\begin{align}
\nonumber
\EE(M_i-M_{i+1})^2
&=
\EE\left[
V_i^{-1/\alpha}\left(
1-\left(\frac{V_{i+1}}{V_i}\right)^{-1/\alpha}
\right)
\right]^2
\\
\nonumber
&=
\EE\left[
V_i^{-2/\alpha}\left(
1-\left(1+\frac{W_{i+1}}{V_i}\right)^{-1/\alpha}
\right)^2
\right]
\\
\nonumber
&\leq
\EE\left[
\frac{1}{\alpha^2}V_i^{-2/\alpha}V_i^{-2}W_{i+1}^2
\right]
\\
\nonumber
&=C\EE V_i^{-2-2/\alpha},
\end{align}
for some constant $C$, since $V_i$ and $W_{i+1}$ are independent.
Arguing as at (\ref{Gammacalc}) and (\ref{Gammafact}),
we have that $V_i$ has Gamma$(i,1)$ distribution, and we 
obtain that for some constant $\tC$ and all large enough $i$ 
\begin{equation}
\label{Mbound}
\EE(M_i-M_{i+1})^2\leq \tC i^{-2-2/\alpha}.
\end{equation}

Suppose $k$ is large enough that (\ref{Mbound}) holds for all $i>k$.
Then using Cauchy-Schwarz, we obtain that for all $\epsilon>0$,
\begin{align*}
U_k
&=\sum_{i=k+1}^\infty L_i(M_i-M_{i+1})
\\
&\leq \left[
\sum_{i=k+1}^\infty
\left(
i^{-1/2-\epsilon}
\right)^2
\right]^{1/2}
\left[
\sum_{i=k+1}^\infty
\left(
i^{1/2+\epsilon}
L_i
(M_i-M_{i+1})
\right)^2
\right]^{1/2}.
\end{align*}
The first sum is finite for any $\epsilon>0$.
Thus squaring and taking expectations, we have that 
\begin{align*}
\EE U_k^2
&\leq
c'\sum_{i=k+1}^\infty
i^{1+\epsilon}\EE L_i^2 \EE(M_i-M_{i+1})^2
\\
&\leq
c''\sum_{i=k+1}^\infty i^{1+\epsilon} i i^{-2-2/\alpha}
\\
&=c''\sum_{i=k+1}^\infty i^{-2/\alpha+\epsilon},
\end{align*}
for some constants $c', c''$. Since $\alpha<2$,
this sum is finite for small enough $\epsilon$,
and so $\EE U_k^2<\infty$ as desired.\enpf

\begin{proposition}
\label{prop:holder}
For all $\beta<\alpha$ and all $\bfu, \bfv\in\RR_+^2$,
we have
\[ 
\be |T(\bfu)-T(\bfv)|^{\beta} \leq
4\max\{\|\bfu\|,\|\bfv\|\}^{\beta/\alpha} \|\bfu-\bfv\|^{\beta/\alpha}
\be T({\bf 1})^{\beta}. 
\]
\end{proposition}

\noindent
\emph{Proof:} Recall that $T(\bfu)$ is increasing in the partial order
on $\br_+^2$. We just need to consider the two cases where $\bfu \leq
\bfv$ and where they are not comparable, so that, say, $u_1<v_1$ and
$u_2>v_2$. 

For $\bfu\leq\bfv$, it is a simple observation that 
\begin{equation}
\label{tvtubound}
T(\bfv) -T(\bfu) \leq T\big((0,u_2),\bfv\big) + T\big((u_1,0),\bfv\big), 
\end{equation}
where, for $\bx\leq\by\in\RR_+^2$,
$T(\bx,\by)$ denotes the maximal weight of an increasing path 
from $\bx$ to $\by$.
By the scaling in the field we have the distributional relationships
\begin{equation}\label{scalingrel}
T\big((a,b),(c,d)\big) 
\stackrel{d}{=} 
T(c-a,d-b) 
\stackrel{d}{=} 
(c-a)^{1/\alpha}(d-b)^{1/\alpha}T(1,1).
\end{equation}
For any random variables $X_1$, $X_2$, we have
$\EE(X_1+X_2)^\beta\leq 2^{\max\{\beta,1\}}
\max\{\EE X_1^\beta, \EE X_2^\beta\}$.
Thus from (\ref{tvtubound}) and (\ref{scalingrel}) we obtain
\begin{align}
\nonumber
\be |T(\bfv)-T(\bfu)|^{\beta} 
&\leq 2^{\max\{\beta,1\}}
\max\left\{
\EE T(v_1-u_1,v_2)^\beta,
\EE T(v_1,v_2-u_2)^\beta
\right\}
\\
\label{firstform}
&\leq
4 \EE T({\bf 1})^\beta \max\left\{
(v_1-u_1)^{\beta/\alpha}v_2^{\beta/\alpha},
v_1^{\beta/\alpha}(v_2-u_2)^{\beta/\alpha}
\right\}.
\end{align}
For the case $u_1<v_1$, $v_2<u_2$, we observe similarly that
\begin{equation}
\label{eq:tuv3}
|T(\bv)-T(\bu)|\leq T\big((u_1,0), (v_1,u_2)\big)+T\big((0,v_2),(v_1,u_2)\big).
\end{equation}
Raising to the power $\beta$ and taking expectations as above,
we obtain that 
\begin{equation}\label{secondform}
\EE |T(\bv)-T(\bu)|^\beta
\leq 4\EE T({\bf 1})^{\beta}\max\left\{
(v_1-u_1)^{\beta/\alpha}v_2^{\beta/\alpha},
(u_2-v_2)^{\beta/\alpha}u_1^{\beta/\alpha}
\right\}.
\end{equation}
Combining the estimates from 
(\ref{firstform}) and (\ref{secondform}) now gives the result.
\enpf

We are now ready to define our analogue of the Airy process. If we set
$H_u = \Theta(u,-u) = T(e^u,e^{-u})$, we have a one-dimensional stationary
process $\{H_u,u\in\br\}$, 
as the processes $(H_{u+t}, u\in\br)$ and $(H_u, u\in\br)$ 
have the same distribution for all
$t\in \br$. We note that as $H_0 = T(1,1)$, the marginal
distribution for our stationary process is the distribution of the
limit random variable in our continuous model for heavy-tailed last
passage percolation. By applying the estimates for the random field we
have estimates on the H\"older continuity of the heavy-tailed Airy
process. 

\begin{corollary}
For $\beta<\alpha$, we have: 
\begin{itemize}
\item[(i)] $\be H_0^\beta <\infty$.
\item[(ii)] For each $\tau>0$ and all $u,v\in[-\tau,\tau]$,  
\[ \be |H_u-H_v|^{\beta} \leq 2e^{2\beta\tau/\alpha}
|u-v|^{\beta/\alpha} \be H_0^\beta. \]
\end{itemize}
\end{corollary}

\noindent
\emph{Proof:} The first part is
Proposition~\ref{prop:fieldmoment}. The second follows from the second
part of the proof of Proposition~\ref{prop:holder} as, assuming $u<v$
and using (\ref{eq:tuv3}) with $v_1=e^v, v_2=e^{-v}, u_1=e^u,u_2=e^{-u}$,
\begin{eqnarray*}
 \be |H_u-H_v|^{\beta} &\leq& (e^v(e^{-u}-e^{-v}))^{\beta/\alpha} \be
H_0^{\beta} + ((e^v-e^u)e^{-u})^{\beta/\alpha} \be H_0^{\beta} \\
&=& 2 (e^{v-u}-1)^{\beta/\alpha} \be H_0^{\beta} \\
&\leq & 2(v-u)^{\beta/\alpha} e^{(v-u)\beta/\alpha}\be H_0^{\beta},
\end{eqnarray*}
and the result follows. \enpf

Finally we give a weak convergence result for our 
heavy-tailed Airy process.

\begin{theorem}
The sequence $\{a_{n^2}^{-1} T^{(n)}(e^u,e^{-u})\}_{u\in [-\tau,\tau]}$
converges weakly to $\{H_u\}_{u\in [-\tau,\tau]}$ in $D[-\tau,\tau]$.
\end{theorem}

\noindent
\emph{Proof:}
We follow the approach outlined in \cite{Resnick86} for such
weak convergence problems, in particular the proof of
\cite{Resnick86}~Proposition~3.4.

Since we restrict to $u\in[-\tau, \tau]$,
we can consider only those points of the 
PRM whose locations fall in the box $[0,c_x]\times[0,c_y]$
where $c_x=c_y=e^\tau$. Thus we can write
the points as a sequence $(Y_i, Z_i), i=1,2,\dots$
such that the sequence $(Z_i)$ is decreasing and such that
$Y_i\in [0,c_x]\times [0,c_y]$ for all $i$.

Define $H^{k}_u = T_k(e^u,e^{-u})$ where 
\[ T_k(x,y) = \sup_{A\in \cC^k_{xy}} \sum_{i\in A} Z_i, \]
with 
\[ \cC^k_{xy} = \{A \subset\{1,\dots,k\} : Y_i\sim Y_j \,\forall
i,j\in A\}. \]

As before, we need a corresponding formulation of the
discrete model. For given $n$,
we work on the box $B^n=\{1,\dots,\lceil ne^\tau\rceil\}^2$.
The weight at a point $y\in B^n$ has distribution $F$,
independently for different points. 
We represent the weights and their positions
by a vector $(Y^n_i, M^n_i,1\leq i\leq \lceil ne^\tau\rceil^2)$,
where $M^n_i$ form a decreasing sequence. 
The interpretation is that $Y^n_i$ is the location
of the $i$th largest weight $M^n_i$ in the box $B^n$.

We can then define 
$H^{(n),k} = 
\tilde{T}^{(n)}_k(\lceil ne^u\rceil ,\lceil ne^{-u}\rceil)$ 
where 
\[ \tilde{T}^{(n)}_k(\lceil nx \rceil ,\lceil ny \rceil) 
= \sup_{A\in \cC^{(n),k}_{xy}} \sum_{i\in A}
a_{n^2}^{-1} M^{(n)}_i, \]
with 
\[ \cC^{(n),k}_{xy} = \{A \subset\{1,\dots,k\wedge n^2\} : 
\forall i,j\in A, Y^{(n)}_i \sim Y^{(n)}_j\}. \]

For the proof we will establish that $H^{(n),k}$ converges weakly
to $H^k$, that $H^k \to H$ locally uniformly and that 
\[ \lim_{k\to\infty} \limsup_{n\to\infty}
\bp(\rho(H^{(n),k},H^{(n)})>\epsilon) = 0, \]
for each $\epsilon>0$, where $\rho$ is the metric for the Skorohod
topology. 


We begin by establishing $H^k \to H$ as $k\to\infty$. This is
a consequence of the construction via a PRM. For each $u$, $H^k_u$ is an
increasing function of $k$ and converges to $H_u$. 
With probability 1, this holds 
uniformly for each $u\in[-\tau, \tau]$, 
since for all such $u$ we have 
\[
0\leq H_u-H_u^k \leq 
\sup_{A\in \cC_{e^\tau, e^\tau}}\sum_{i\in A, i>k} M_i;
\]
the upper bound is finite with probability 1 exactly as in 
Lemma \ref{Slemma}.

Next we wish to show the 
weak convergence of $H^{(n),k}$ to $H^k$.
This can be done by an extension of the method 
of Proposition \ref{coupledprop}.
Note that with probability 1, 
no two points of the PRM share a vertical coordinate 
or a horizontal coordinate, and in addition
no point of the PRM falls on the line 
parametrised by $(x,y)=(e^u, e^{-u})$.
Then under the same couplings used in the proof 
of Proposition \ref{coupledprop},
one obtains that with probability 1,
$H^{(n),k}\to H^k$ in the Skorohod space. 
Thus this weak convergence also holds as desired.

Finally we need to control $\rho(H^{(n),k},H^{(n)})$. 
Fix an $\epsilon>0$.
Consider
\begin{eqnarray*}
\bp\left(\rho(H^{(n),k},H^{(n)}) > \epsilon\right) 
&\leq & 
\bp\left(\sup_{-\tau<u<\tau}
\left|H^{(n),k}_u - H^{(n)}_u\right|>\epsilon\right) 
\\
&=&  
\bp\left(\sup_{-\tau<u<\tau} \left|\tilde{T}_k^{(n)}(e^u,e^{-u}) -
\tilde{T}^{(n)}(e^u,e^{-u})\right|>\epsilon\right) 
\\ 
&=&  
\bp\left(\sup_{-\tau<u<\tau} 
\left|
\sup_{A\in C^{(n),k}_{e^u,e^{-u}}} \sum_{i\in A}
a_{n^2}^{-1} M_i^{n} - \sup_{A\in C^{(n)}_{e^u,e^{-u}}} \sum_{i\in A}
a_{n^2}^{-1} M_i^{n} 
\right|
>\epsilon\right) 
\\
&\leq&  
\bp\left(\sup_{-\tau<u<\tau} 
\left| 
\tilde{S}^{(n)}_k(e^u,e^{-u})
\right|
>\epsilon\right),
\end{eqnarray*}
where $\tilde{S}^{(n)}_k(e^u,e^{-u}) = \sup_{A\in
C^{(n)}_{e^u,e^{-u}}} \sum_{i\in A, i>k} a_{n^2}^{-1} M_i^{n}$. By
monotonicity we have
\[  \sup_{-\tau<u<\tau} \left| \tilde{S}^{(n)}_k(e^u,e^{-u})\right| \leq
\tilde{S}^{(n)}_k(e^\tau,e^\tau), \]
and using the scaling
\[ \bp\Big(\rho(H^{(n),k},H^{(n)}\Big) > \epsilon) \leq
\bp\Big(\tilde{S}^{(n)}_{k}(1,1) > \epsilon e^{-2\tau/\alpha}\Big). \]

By Proposition~\ref{Aprop} we have that
$\bp\big(\tilde{S}^{(n)}_k > \epsilon e^{-2\tau/\alpha}\big) \to 0$ 
as $k\to \infty$ uniformly in $n$. 
Thus indeed we have
\[ \lim_{k \to\infty} \limsup_{n\to \infty}
  P\big(\rho(H^{(n),k},H^{(n)}\big) >\epsilon) =0. \]
Putting the three pieces together we have shown the weak convergence. \enpf

\section{Higher-dimensional heavy-tailed last passage percolation}
\label{dsection}
Up to this point we have considered only two-dimensional models.
In this section we indicate how to extend most of the results 
to higher dimensions in a natural way.

For general $d\geq 2$, we consider the passage time
from the point $(1,1,\dots,1)$ to the point $(n,n,\dots, n)$.

We now consider a sequence of locations $Y_i^{(n)}$,
$i=1,2,\dots, n^d$ which form a uniform random
permutation of the set $\{1/n, 2/n, \dots, 1\}^d\subset [0,1]^d$,
and a corresponding sequence of weights $M_i^{(n)}$, $i=1,2,\dots,n^d$
which are given by the order statistics, in decreasing order,
of a sample of size $n^d$ from the underlying weight distribution $F$.
We now assume that the tail of $F$ is regularly varying
with index $\alpha<d$.

Defining $\cC^{(n)}$ as before, we set 
$T^{(n)}=\sup_{A\in\cC^{(n)}} \sum_{i\in A} M_i^{(n)}$,
and $\tT^{(n)}=a_{n^d}^{-1} T^{(n)}$.

The continuous model is defined just as before;
the locations $Y_i$ are now drawn i.i.d.\ and uniformly at random
from the box $[0,1]^d$ rather than the square $[0,1]^2$.

Then $\tT^{(n)}\to T$ in distribution as $n\to\infty$;
the method of proof is essentially identical to that
used for Theorem \ref{Ttheorem} in the case $d=2$.

The multivariate extensions described in Sections 
\ref{preAiry} and \ref{Airysection}
go through in an analogous way. 
For example, we can now obtain a process $\Theta$ which is 
stationary on $\RR^d$ such that 
\[ 
\Bigg\{ \exp\left(-\frac{u_1+\dots+u_d}{\alpha}\right) 
a_{n^d}^{-1}T^{(n)}\left(e^{u_1},\dots,e^{u_d}\right), \bu\in\RR^d\Bigg\}
\to \Big\{\Theta(u_1,\dots,u_d), \bu\in\RR^d\Big\}
\]
as $n\to\infty$, in the sense of convergence of
finite-dimensional distributions;
here $T^{(n)}(u_1,\dots,u_d)$ is the 
maximal weight of a path from $(1,\dots,1)$  to 
the point $(\lceil n u_1 \rceil,
\dots, \lceil n u_d \rceil)$.

We turn to the path convergence as developed in Section
\ref{pathsection}. 
Proposition \ref{uniquenessprop} and Theorem \ref{paththeorem}
extend easily, with the same method of proof. 
However, extending Theorem \ref{paththeorem2}, 
concerning the convergence of optimal paths 
viewed as random subsets of $[0,1]^d$,
is more problematic. 
Again we are unable to prove that the optimal path 
for the continuous model (i.e.\ the closure of 
$\bigcup_{i\in A^*} Y_i$) is connected
(although we expect this to be true).
In the case $d=2$ this caused a little inconvenience 
but we could work around it by observing that
any ``jumps'' in the path consist of 
horizontal or vertical line segments, and hence that
at least there exists a unique connected increasing path
that contains the optimal path.

For $d\geq 3$, however, a jump could, for example, 
cross a square of zero volume in $\RR^d$ but 
with non-zero area. Then there is no longer a unique 
way to extend the optimal path to a connected increasing path,
and thus there is an ambiguity in the limit object.
If we could prove the conjecture that the optimal path 
itself is connected, the convergence in distribution 
of the discrete optimal paths $P^{(n)*}$ would follow as before.

In the case $\alpha=0$, we can in fact prove the connectedness
of the optimal path for the continuous model 
(defined as in Section \ref{greedysection} using
the ``greedy algorithm''). 
It is not clear how to extend the multifractal 
analysis of Section \ref{greedysubsection}.
However, by analysing a branching
random walk associated with the algorithm which constructs
the greedy path, one can obtain
that the function from, say, $x_1\in[0,1]$
to $(x_2,\dots, x_d)\in[0,1]^{d-1}$ which describes the path
is almost surely everywhere continuous and strictly increasing
(although a.s.\ it also has derivative 0 almost everywhere).
Thus one can show that $P^{(n)*}\to P^*$ 
in distribution (under the Hausdorff metric on subsets of $[0,1]^d$)
for all $d$ in the case $\alpha=0$.

\section*{Acknowledgments}
We are grateful for the support of the Isaac Newton Institute in Cambridge; 
this work began during the programme \textit{Interaction
and Growth in Complex Stochastic Systems}.


\bigskip

\parbox{0.52\linewidth}
{
\textsc{Mathematical Institute,\\
University of Oxford,\\
24-29 St Giles,\\
Oxford OX1 3LB,\\
UK}\\
\texttt{hambly@maths.ox.ac.uk}\\
\texttt{http://www.maths.ox.ac.uk/$\tilde{\,\,\,\,}$hambly}}
\hfill
\parbox{0.47\linewidth}
{\textsc{Department of Statistics,\\
University of Oxford,\\
1 South Parks Road,\\
Oxford OX1 3TG\\
UK}\\
\texttt{martin@stats.ox.ac.uk}\\
\texttt{http://www.stats.ox.ac.uk/$\tilde{\,\,\,\,}$martin}}

\end{document}